\documentclass{article}

\usepackage{amsfonts}
\usepackage{color}
\usepackage{graphics}
\usepackage{epsfig,psfrag,graphicx}
\usepackage{amssymb}
\usepackage{eepic,epic}
\usepackage{epsfig} 
\usepackage{amsmath} 
\usepackage{amssymb} 
\usepackage{amsbsy} 
\usepackage{fullpage}

\def\Xint#1{\mathchoice
{\XXint\displaystyle\textstyle{#1}}%
{\XXint\textstyle\scriptstyle{#1}}%
{\XXint\scriptstyle\scriptscriptstyle{#1}}%
{\XXint\scriptscriptstyle%
\scriptscriptstyle{#1}}%
\!\int}
\def\XXint#1#2#3{{\setbox0=\hbox{$#1{#2#3}{%
\int}$ }
\vcenter{\hbox{$#2#3$ }}\kern-.6\wd0}}

\def\dashint{\Xint-}

\newtheorem{Theorem}{Theorem}[section]
\newtheorem{Proposition}{Proposition}[section]
\newtheorem{Lemma}{Lemma}[section]
\newtheorem{Corollary}{Corollary}[section]

\newtheorem{Remark}{Remark}[section]

\newcommand{\Div}       {{\rm div}_x}
\newcommand{\dive}       {{\rm div} \:}

\newcommand{\pa} {\partial}

\newcommand{\dx}       		 {{\rm \,d}x}

\newcommand{\dy}			{{{\rm \,d}y}}

\newcommand{\ep}    {\varepsilon}
\newcommand{\eps}    {\varepsilon}

\newcommand{\R}         {\mathbb{R}}
\newcommand{\N}         {\mathbb{N}}
\newcommand{\T}    {{\mathbb T}}
\newcommand{\Z}         {\mathbb{Z}}

\newcommand{\tS}		{{{S}}}

\newcommand{\na}{\nabla}

\newcommand{\bS}   {{S}}
\newcommand{\tA}		{{{A}}}
\newcommand{\tD}		{{{D}}}

\newcommand{\ub}		{{{u}_{bl}}}

\newcommand{\we}		{{{w}^\ep}}

\newcommand{\ue}        {{{u}^\ep}}
\newcommand{\uz}		{{{u}^0}}

\newcommand{\qed}{\hfill $\Box $ \bigskip }

\def\bbbone{{\mathchoice {\rm 1\mskip-4mu l}
{\rm 1\mskip-4mu l} {\rm 1\mskip-4.5mu l} {\rm 1\mskip-5mu l}}}

\date{\today
}

\begin{document}

\title{Boundary layer for a  non-Newtonian  flow \\ over a rough surface}

\author{David G\'erard-Varet
\thanks{Institut de Math\'{e}matiques de Jussieu et Universit\'{e} Paris 7, 175 rue du Chevaleret, 75013 Paris France
({\tt gerard-varet@math.jussieu.fr})}
, Aneta Wr\'oblewska-Kami\'nska 
\thanks{ Institute of Mathematics, Polish Academy of Sciences, ul. \'Sniadeckich 8, 00-956 Warszawa, Poland
({\tt awrob@impan.pl})}
}

\maketitle

\bigskip

\section{Introduction}
The general concern of this paper is the effect of rough walls on fluids. This effect is important at various scales. For instance, in the area of microfluidics, recent experimental works have emphasized the role of hydrophobic rough walls in the improvement of slipping properties of microchannels. Also, in geophysics, as far as large scale motions are concerned, topography or shore variations can be assimilated to roughness. 
For high Reynolds number flows, an important issue is to understand how localized roughness triggers  instabilities, and transition to turbulence. For laminar flows, the point is rather to understand how distributed  roughness may have  a macroscopic impact on the dynamics. More precisely,  the hope is to be able to encode an averaged effect through an effective boundary condition at a smoothened wall. Such boundary condition, called {\em a wall law}, will avoid to simulate the small-scale dynamics that takes place in a boundary layer in the vicinity of  the rough surface. 

\medskip
The derivation of wall laws for laminar Newtonian flows has been much studied, since the pioneering works of Achdou, Pironneau and Valentin \cite{Achdou:1995, Achdou:1998}, or J\"ager and Mikeli\'c \cite{Mikelic2001,Jager:2003}. See also \cite{Luchini:1995,Amirat:2001a,GV2003,Bresch,Mikelic2013}.  A natural mathematical approach of this problem is by homogenization techniques, the roughness being modeled by a small amplitude/high frequency oscillation. Typically, one considers a Navier-Stokes flow in a channel $\Omega^\eps$ with a rough bottom:  
$$\Omega^\ep = \Omega \cup \Sigma_0 \cup R^\ep. $$
Precisely:  
\begin{itemize}
\item $\Omega = (0,1)^2$ is the flat portion of the channel.  
\item $R^\eps$ is the rough portion of the channel: it reads 
$$ R^\eps = \{ x = (x_1,x_2), x_1 \in (0,1), 0 > x_2 > \eps \gamma(x_1/\eps) \}$$
with a bottom surface $\Gamma^\eps := \{x_2 = \eps \gamma(x_1/\eps) \}$ parametrized by $\eps \ll 1$. Function $\gamma = \gamma(y_1)$ is the {\em roughness pattern}.
\item Eventually, $\Sigma_0 := (0,1) \times \{0\}$ is the interface between the rough and flat part. It is the artificial boundary at which the wall law is set. 
\end{itemize}
Of course, within such model, the goal is to understand the asymptotic behavior of the Navier-Stokes solution $u^\eps$ as $\eps \rightarrow 0$. Therefore, the starting point is a formal approximation  of $u^\eps$ under the form 
\begin{equation} \label{blexpansion} 
 u^\eps_{app}(x)  = u^0(x) + \eps u^1(x)  + \dots  +  u^0_{bl}(x/\eps) + \eps u^1_{bl}(x/\eps) + \dots  .
 \end{equation}  
 In this expansion, the $u^i = u^i(x)$ describe the large-scale part of the flow, whereas the $u^i_{bl} = u^i_{bl}(y)$ describe the boundary layer. The  typical variable $y=x/\eps$ matches the small-scale variations induced by the roughness. In the case of homogeneous Dirichlet conditions at $\Gamma^\eps$, one can check formally that: 
 \begin{itemize}
\item $u^0$ is the solution of the Navier-Stokes equation in $\Omega$, with Dirichlet condition at $\Sigma_0$. \item  $u^0_{bl} = 0$, whereas $u^1_{bl}$ satisfies a Stokes equation in variable $y$ in the boundary layer domain 
 $$\Omega_{bl} := \{y = (y_1, y_2), y_1 \in \R, y_2 > \gamma(y_1)\}.$$  
 \end{itemize}
The next step is to solve this boundary layer system, and show convergence of $u^1_{bl}$ as $y_2 \rightarrow +\infty$ to a constant field $u^\infty = (U^\infty, 0)$. This in turn determines the appropriate boundary condition for the large scale correction $u^1$. From there, considering the large scale part $u^0 + \eps u^1$, one can show that:
\begin{itemize}
\item The limit wall law is a homogeneous Dirichlet condition. Let us point out that this feature persists even starting from a microscopic pure slip condition, under some non-degeneracy of the roughness: \cite{Bucur:2008,Bucur:2012,BoDaGe}.
\item The $O(\eps)$ correction to this wall law is a slip condition of Navier type, with $O(\eps)$ slip length. 
\end{itemize}
All these steps were completed in aforementioned articles, in the case of periodic roughness pattern $\gamma$ : $\gamma(y_1 + 1) = \gamma(y_1)$. Over the last years, the first author has extended this analysis to general patterns of roughness, with ergodicity properties (random stationary distribution of roughness, {\it etc}). We refer to \cite{BaGeVa2008, DGV:2008,GeVaMa2010}. See also \cite{DaGeVa2011}  for some recent work on the same topic.  
 
 \medskip
 The purpose of the present paper is to extend the former analysis to non-Newtonian flows. This  may have various sources of interest. One can think of engineering applications, for instance lubricants to which  polymeric additives confer a shear thinning behavior. One can also think of glaciology: as the interaction of glaciers with the underlying rocks is unavailable, wall laws can help. From a mathematical point of view, such examples may be described by a power-law model. Hence, we consider a system of the following form: 
 \begin{equation} \label{EQ1}
 \left\{ 
 \begin{aligned}
-\dive S(Du)  + \na p  = e_1 & \quad \mbox{in} \: \Omega^\eps, \\
\dive u  = 0 &   \quad \mbox{in} \: \Omega^\eps, \\
u\vert_{\Gamma^\eps} = 0, \quad u\vert_{x_2 = 1}  = 0&, \quad u \: \mbox{$1$-periodic in $x_1$}. 
\end{aligned}
\right. 
\end{equation} 
As usual, $u = u(x) \in \R^2$ is the velocity field, $p = p(x) \in \R$ is the pressure. The source term $e_1$ at the right-hand side of the first equation corresponds to a constant pressure gradient $e_1 = (1,0)^t$ throughout the channel. Eventually, the left-hand side involves the stress tensor of the fluid. As mentioned above, it is taken of power-law type:  $S : \R^{2\times 2}_{\rm sym} \to \R^{2\times 2}_{\rm sym}$ is given by 
\begin{equation} \label{defstress}
S :  	 \R^{2\times 2}_{\rm sym} \to \R^{2\times 2}_{\rm sym}, \quad S(A) = \nu |A|^{p-2}A, \quad \nu > 0, \quad 1 < p < +\infty ,
\end{equation}
where $|A| = (\sum_{i,j} a_{i,j}^2)^{1/2}$ is the usual euclidean norm of the matrix $A$. {\em For simplicity, we shall take $\nu = 1$}. Hence, $S(Du) =  |Du|^{p-2} Du$, where we recall that $Du = \frac{1}{2} (\na u + (\na u)^t)$ is the symmetric part of the jacobian. Following classical terminology, the case $p < 2$ resp. $p > 2$ corresponds to {\em shear thinning} fluids, resp. {\em shear thickening} fluids. The limit case $p=2$ describes a Newtonian flow. Note that we complete the equation in system \eqref{EQ1} by a standard no-slip condition at the top and bottom boundary of the channel. For the sake of simplicity, we assume periodicity in the large scale horizontal variable $x_1$. Finally, we also make a simplifying periodicity assumption on the roughness pattern $\gamma$: 
\begin{equation} 
\mbox{$\gamma$ is $C^{2,\alpha}$ for some $\alpha  >0$, has values in $(-1,0)$, and is  $1$-periodic in $y_1$}.
\end{equation}
For every $\eps > 0$ and any value of $p$, the generalized Stokes system \eqref{EQ1} has a unique solution 
$$ (u^\eps, p^\eps) \in W^{1,p}(\Omega^\ep) \times L^{p'}(\Omega^\ep)/\R .$$
The main point is to know about the asymptotic behavior of $u^\eps$, precisely to build some good approximate solution. With regards to the Newtonian case, we anticipate that this approximation will take a form close to  \eqref{blexpansion}.  Our plan is then: 
\begin{itemize}
\item to derive the equations satisfied by the first terms of expansion \eqref{blexpansion}. 
\item to solve these equations, and show convergence of the boundary layer term to a constant field away from the boundary.
\item to obtain error estimates for the difference between $u^\eps$ and $u^\eps_{app}$. 
\item to derive from there appropriate wall laws.
\end{itemize}
This program will be more difficult to achieve for non-Newtonian fluids, in particular for the shear thinning case $p < 2$, notably as regards the study of the boundary layer equations on $u_{bl} \: :=  u^1_{bl}$. In short, these equations will be seen to read 
$$ - \dive(S(A + D u_{bl})) + \na p = 0, \dive u = 0, \quad y \in \Omega_{bl} $$
for some explicit matrix $A$, together with periodicity condition in $y_1$ and a homogeneous Dirichlet condition at the bottom of $\Omega_{bl}$. Due to the nonlinearity of these equations and the fact that $A \neq 0$, the analysis will be much more difficult than in the Newtonian case, notably the proof of the so-called Saint-Venant estimates. We refer to section \ref{subsecstvenant} for all details. 

\medskip
Let us conclude this introduction by giving some references on related problems. In \cite{MarusicPaloka2000}; E. Maru\v{s}i\'c-Paloka  considers power-law fluids with convective terms in infinite channels and pipes (the non-Newtonian analogue of the celebrated Leray's problem). After an appropriate change of unknown, the system studied in \cite{MarusicPaloka2000} bears some strong similarity to our boundary layer system. However, it is different at two levels : first, the analysis is restricted to the case $p > 2$. Second, our lateral  periodicity condition in $y_1$ is replaced by a no-slip condition. This allows to use Poincar\'e's inequality in the transverse variable, and control  zero order terms (in velocity $u$) by $\na u$, and then by $D u$ through the Korn inequality. It simplifies in this way the derivation of exponential convergence of the boundary layer solution (Saint-Venant estimates). The same simplification holds in the context of paper \cite{BoGiMa-Pa}, where  the behaviour of a Carreau flow through a thin filter is analysed. The corrector describing the behaviour of the fluid near the filter is governed by a kind of boundary layer type system, in a slab that is infinite vertically in both directions. In this setting, one has $A = 0$, and the authors refer to \cite{MarusicPaloka2000} for well-posedness and qualitative behaviour. We also refer to the recent article \cite{Suarez-Grau_2015} dedicated  to power-law fluids in thin domains, with Navier condition and anisotropic roughness (with a wavelength that is  larger than the amplitude). In this setting, no boundary layer analysis is needed, and the author succeeds to describe the limit asymptotics by the unfolding method. Finally, we point out the very recent paper  \cite{ChupinMartin2015}, where an Oldroyd fluid is considered in a rough channel. In this setting, no nonlinearity is associated to the boundary layer, which  satisfies a Stokes problem.

\section{Boundary layer analysis}
From the Newtonian case, we expect the solution $(u^\eps, p^\eps)$ of  \eqref{EQ1} to be approximated by  
$$ u^\eps \approx u^0(x) + \eps u_{bl}(x/\eps), \quad p^\eps \approx p^0(x) +  p_{bl}(x/\eps) ,$$
where
\begin{itemize} 
\item $(u^0, p^0)$ describes the flow away from the boundary layer. We shall take $u^0 = 0$ and $p^0 = 0$ in the rough part $R^\eps$ of the channel.
\item $(u_{bl}, p_{bl}) = (u_{bl}, p_{bl})(y)$ is a boundary layer corrector defined on the slab  
 $$\Omega_{bl} := \{y = (y_1, y_2), y_1 \in \T, y_2 > \gamma(y_1)\},$$    
where $\T$ is the torus $\R/\Z$. This torus corresponds implicitly to a periodic boundary condition in $y_1$, which is inherited from the periodicity of the roughness pattern $\gamma$.  
We denote 
$$ \Omega_{bl}^{\pm} \: := \: \Omega_{bl} \cap \{ \pm y_2 > 0 \} $$
its upper and lower parts, and 
$$ \Gamma_{bl} := \{ y = (y_1, y_2), y_1 \in \T, y_2 = \gamma(y_1)\} $$
its bottom boundary.  As the boundary layer corrector is supposed to be localized, we expect that
$$ \na u_{bl} \rightarrow 0 \quad \mbox{as $y_2 \rightarrow +\infty$}. $$
\end{itemize}
With this constraint in mind, we take $(u^0, p^0)$ to be the solution of  
 \begin{equation} \label{EQ2}
 \left\{ 
 \begin{aligned}
-\dive S(Du^0)  + \na p^0  = e_1 & \quad \mbox{in} \: \Omega, \\
\dive u^0  = 0 &   \quad \mbox{in} \: \Omega, \\
u^0\vert_{\Sigma_0} = 0, \quad u^0\vert_{x_2 = 1}  = 0&, \quad u^0 \mbox{ $1$-periodic in $x_1$}. 
\end{aligned}
\right. 
\end{equation} 
The solution is explicit and generalizes the Poiseuille flow. A simple calculation yields: for all $x \in \Omega$,  
\begin{equation}  \label{Poiseuille}
p^0(x) = 0, \quad u^0(x) = (U(x_2), 0), \quad U(x_2) = \frac{p-1}{p} \left(\sqrt{2}^{-\frac{p}{(p-1)}} -  \sqrt{2}^{\frac{p}{(p-1)}} \left|x_2 - \frac{1}{2}\right|^{\frac{p}{p-1}} \right). 
\end{equation}
We extend this solution to the whole rough channel by taking: $u^0 = 0, p^0 = 0$ in $R^\eps$. 
This zero order approximation is clearly continuous across the interface $\Sigma_0$, but the associated stress  is not: denoting 
\begin{equation} \label{defA} 
A \: := D(u^0)\vert_{y_2 = 0^+}\: = \frac{1}{2}\begin{pmatrix} 0 & U'(0) \\ U'(0) & 0 \end{pmatrix}, \quad \mbox{with $U'(0) = \sqrt{2}^{\frac{p-2}{p-1}}$} 
\end{equation} 
we obtain 
$$ [ S(Du^0) n - p^0 n ]\vert_{\Sigma_0} = |A|^{p-2}A n = \left( \begin{smallmatrix} - \sqrt{2}^{-p} U'(0)^{p-1}  \\ 0 \end{smallmatrix}\right)  =  \left( \begin{smallmatrix} - \frac{1}{2}  \\ 0 \end{smallmatrix}\right)  $$ 
with $n = - e_2 = -(0,1)^t$ and $[f] := f\vert_{x_2 = 0^+} - f\vert_{x_2 = 0^-}$. 
 
 \medskip
 This  jump should be corrected by $u_{bl}$, so that the total approximation $u^0(x) + \eps u_{bl}(x/\eps)$ has no jump. This explains the amplitude $\eps$ of the boundary layer term, as its gradient will then be $O(1)$. By Taylor expansion $U(x_2) = U(\eps y_2) = U(0) + \eps U'(0) y_2 + \dots$  we get formally $D(u^0 + \eps u_{bl}(\cdot/\eps)) \approx A + D u_{bl}$, where the last symmetric gradient is with respect to variable $y$. We then derive the following boundary layer system: 
 \begin{equation} \label{BL1}
 \left\{
 \begin{aligned}
 - \dive S(A + D u_{bl}) + \na p_{bl}  & = 0 \quad \mbox{in} \: \Omega_{bl}^+, \\
   - \dive S(D u_{bl}) + \na p_{bl}  & = 0 \quad \mbox{in} \: \Omega_{bl}^-, \\
 \dive u_{bl} & = 0 \quad \mbox{in} \: \Omega_{bl}^+ \cup \Omega_{bl}^-, \\
u_{bl}\vert_{\Gamma_{bl}} & = 0,  \\
u_{bl}\vert_{y_2 = 0^+} - u_{bl}\vert_{y_2 = 0^-}  & = 0,   
 \end{aligned}
 \right.
 \end{equation}
 together with the jump condition 
\begin{equation} \label{BL2}
	\left( S( A + D u_{bl}) n - p_{bl} n \right) |_{y_2 = 0^+} 
	- \left(S(D u_{bl}) n  - p_{bl} n\right) |_{y_2= 0^-} = 0, \quad n = (0,-1)^t.  
	\end{equation}
Let us recall that the periodic boundary condition in $y_1$ is encoded in the definition of the boundary layer domain. The rest of this section will be devoted to the well-posedness and qualitative properties of \eqref{BL1}-\eqref{BL2}. We shall give detailed proofs only for the more difficult case $p < 2$, and comment briefly on the case $p \ge 2$ at the end of the section. Our main results will be the following:  
\begin{Theorem} {\bf (Well-posedness)} \label{thWP}

\smallskip
\noindent
For all $1 < p < 2$,  \eqref{BL1}-\eqref{BL2} has a unique solution
$ (u_{bl},p_{bl}) \in W^{1,p}_{loc}(\Omega_{bl}) \times L^{p'}_{loc}(\Omega_{bl})/\R  $
satisfying for any $M > |A|$: 
$$D u_{bl} \, 1_{\{|D u_{bl}| \le M\}} \in L^2(\Omega_{bl}), \quad  D u_{bl} \,  1_{\{|D u_{bl}| \ge M\}} \in L^p(\Omega_{bl}).$$

\smallskip
\noindent
For all $p \ge  2$,  \eqref{BL1}-\eqref{BL2} has a unique solution
$ (u_{bl},p_{bl}) \in W^{1,p}_{loc}(\Omega_{bl}) \times L^{p'}_{loc}(\Omega_{bl})/\R  $
s.t.  $D u_{bl} \in L^p(\Omega_{bl}) \cap L^2(\Omega_{bl})$.
\end{Theorem}

 \begin{Theorem} {\bf (Exponential convergence)} \label{thEC}

\smallskip
\noindent
For any $1 < p < +\infty$,  the solution given by the previous theorem converges exponentially, in the sense that for some $C, \delta > 0$ 
$$ | u_{bl}(y) - u^\infty | \le C e^{-\delta y_2} \quad \forall \: y \in \Omega_{bl}^+,$$
 where $u^\infty= (U^\infty, 0)$ is some constant horizontal vector field. 
 \end{Theorem} 
 
 \subsection{Well-posedness} \label{sectionWP}
 
 \subsubsection*{A priori estimates}
 We focus on the case $1 < p < 2$, and provide the {\it a priori} estimates on which the well-posedness is based. The easier case $p \ge 2$ is discussed at the end of the paragraph. As $A$ is a constant matrix, we have from  \eqref{BL1}:  
 $$ - \dive S(A + D u_{bl}) + \dive S(A) + \na p_{bl}  = 0  \quad \mbox{in} \: \Omega_{bl}^+, \quad    - \dive S(D u_{bl}) + \na p_{bl}   = 0 \quad \mbox{in} \: \Omega_{bl}^-. $$
 We multiply the two equations by $D u_{bl}$ and integrate over $\Omega_{bl}^+$ and   
 $\Omega_{bl}^-$ respectively. After integrations by parts, accounting for the jump conditions at $y_2 = 0$, we get 
 \begin{equation} \label{variationalBL}
  \int_{\Omega_{bl}^+} (S(A + Du_{bl}) - S(A)) : D u_{bl} \dy + \int_{\Omega_{bl}^-} S(D u_{bl}) : D u_{bl} \dy = -\int_{y_2 = 0} S(A)n \cdot u_{bl} {\rm\,d}S. 
  \end{equation}
 The right-hand side is controlled using successively Poincar\'e  and Korn inequalities (for the Korn inequality, see the appendix): 
\begin{equation}
  |\int_{y_2 = 0} S(A)n \cdot u_{bl} \dy | \le C \| u_{bl} \|_{L^p(\{ y_2 = 0\})} \le C' \| \na u_{bl} \|_{L^p(\Omega_{bl}^-)} \le C'' \| D u_{bl} \|_{L^p(\Omega_{bl}^-)} .
 \end{equation}
 As regards the left-hand side, we rely on the following vector inequality, established in \cite[p74, eq. (VII)]{Lind}: for all $1 < p \le 2$,  for all vectors $a,b$
\begin{equation}   \label{vectorinequality}
   ( |b|^{p-2}b - |a|^{p-2} a \: | \: b-a) \: \ge \: (p-1) |b-a|^2 \int_0^1 |a + t(b-a)|^{p-2} dt.
\end{equation}  
 In particular, for any $M > 0$,  if $|b-a| \le M$, one has  
 \begin{equation}  \label{ineq1} 
 ( |b|^{p-2}b - |a|^{p-2} a \: | \: b-a) \: \ge \: \frac{p-1}{(|a| + M)^{2-p}} |b-a|^2 ,
 \end{equation} 
 whereas if $|b-a| > M >  |a|$,  we get 
 \begin{equation}  \label{ineq2}
   ( |b|^{p-2}b - |a|^{p-2} a \: | \: b-a) \: 
   \ge (p-1) |b-a|^2 \int_{\frac{|a|}{|b-a|}}^1  \left(2 t |b-a|\right)^{p-2} dt  \ge  2^{p-3} \left(1 - (|a|/M)^{p-1}\right) |b-a|^p.
\end{equation} 
We then apply such inequalities to \eqref{variationalBL}, taking $a = A$, $b = A + Du_{bl}$. For $M > |A|$, there exists $c$ dependent on $M$ such that  
$$
 \int_{\Omega_{bl}^+} (S(A + Du_{bl}) - S(A)) : D u_{bl} \dy \ge c \int_{\Omega_{bl}^+} \bbbone_{\{ |Du_{bl}| \le M \}} |Du_{bl}|^2 \dy  \: + \:   \int_{\Omega_{bl}^+}  \bbbone_{\{ |Du_{bl}| > M \}} |D u_{bl}|^p \dy ,
 $$ 
 so that for some $C$ dependent on $M$
\begin{equation*}
\int_{\Omega_{bl}^+}| \bbbone_{\{ |Du_{bl}| \le M \}} Du_{bl}|^2 \dy  \: + \:   \int_{\Omega_{bl}^+} |\bbbone_{\{ |Du_{bl}| > M \}}  D u_{bl}|^p \dy + \int_{\Omega_{bl}^-}  | D u_{bl}|^p \dy \: \le \: C \, \| D u_{bl} \|_{L^p(\Omega_{bl}^-)} .
\end{equation*}       
Hence, still for some $C$ dependent on $M$:  
\begin{equation} \label{aprioriestimate}
\int_{\Omega_{bl}^+} 1_{\{ |Du_{bl}| \le M \}} |Du_{bl}|^2 \dy  \: + \:   \int_{\Omega_{bl}^+} 1_{\{ |Du_{bl}| > M \}}  |D u_{bl}|^p \dy
+ \int_{\Omega_{bl}^-}  | D u_{bl}|^p \dy \: \le \: C .
\end{equation} 
This is the {\it a priori} estimate on which Theorem  \ref{thWP} can be established (for $p \in ]1,2]$).   Note that this inequality implies  that for any height $h$,
$$ \| Du_{bl} \|_{L^p(\Omega_{bl} \cap \{y_2 \le h \})} \le C_h $$
(bounding the $L^p$ norm by the $L^2$ norm on a bounded set). Combining with  Poincar\'e and Korn inequalities, we obtain that $u_{bl}$ belongs to $W^{1,p}_{loc}(\Omega_{bl})$. 

\medskip 
In the case $p\geq 2$, we can directly use the following inequality,  which holds for all $a,\ b \in \R^n$:
	\begin{equation}\label{abp_1}
	 | a - b |^p 2^{2-p}  \leq  2^{-1}  \left( |b|^{p-2} + |a|^{p-2}\right) |b-a|^{2} \leq\left\langle  |a|^{p-2} a - |b|^{p-2}|b|, a-b \right\rangle .
 	\end{equation}
It provides both an $L^2$ and $L^p$ control of the symmetric gradient of the solution. Indeed, taking $a = \tA + \tD_y \ub$, $b= \tA$  and using \eqref{variationalBL} we get the following a'priori estimates for $p\geq 2$ 
	\begin{equation}\label{apesDup}
	\begin{split}
	&  2^{2-p} \int_{\Omega_{bl}^{+}} | \tD  \ub |^p  \dy +  \int_{\Omega_{bl}^{-}} | \tD \ub |^p  \dy    + 
	 2^{-1} |A|^{p-2}  \int_{\Omega_{bl}^+} | \tD \ub |^2  \dy  \: \\
	 &  {\leq} \:  
	  \int_{\Omega_{bl}^{+}} \left( S(A+ Du_{bl}) - S(A) \right) : D \ub \dy   + \int_{\Omega_{bl}^{-}} S(D u_{bl} ) : D \ub \dy  \\
	& = \: -   \int_{\Sigma_0} S(A) n \cdot u_{bl}  {\rm \, d}S \\
	& \leq \:   c(\alpha)  \| S(A) \|^{p'}_{L^{p'}(\Sigma_0)}  + \alpha \| \ub \|^{p}_{L^p(\Sigma_0)}
	\leq  \: c(\alpha)  \| S(A) \|^{p'}_{L^{p'}(\Sigma_0)}  
	+ \alpha  C_\Gamma \| \nabla \ub \|^{p}_{L^p(\Omega_{bl}^{-})}
	\\ & 
	\leq c(\alpha)  \| S(A) \|^{p'}_{L^{p'}(\Sigma_0)}  
	+ \alpha  C_\Gamma  C_K\| \tD \ub \|^{p}_{L^p(\Omega_{bl}^{-})} ,
	\end{split}
	\end{equation}
where the trace theorem, the Poincar\'e inequality  and the Korn inequality were employed. 
By choosing the coefficient $\alpha$ small enough, and by the imbedding of $L^p(\Omega_{bl}^-)$ in $L^2(\Omega_{bl}^-)$,  \eqref{apesDup} provides 
	\begin{equation}\label{es:Dup}
	\int_{\Omega_{bl}}   |\tD \ub|^{p} +  |\tD \ub|^{2} \dy \leq  C   \| S(A) \|^{p'}_{L^{p'}(\Sigma_0)} < \infty .
	\end{equation}
Eventually,  by  Korn and Poincar\'e inequalities: $\ub \in W^{1,p}(\Omega_{bl})$ for $2\leq p < \infty$.

 \subsubsection*{Construction scheme for the solution}
 We briefly explain how to construct a solution satisfying the above estimates. We restrict to the most difficult case $p \in ]1,2]$. 
There are two steps: 

\medskip
 {\em Step 1}:  we solve the same equations, but   {\em in the bounded domain  $\Omega_{bl,n} = \Omega_{bl} \cap \{ y_2 < n \}$}, with a Dirichlet boundary condition at the top. As 
 $\Omega_{bl,n}$ is bounded, the imbedding of $W^{1,p}$ in $L^p$ is compact, so that a solution $u_{bl,n}$  can be built  in a standard way. Namely,  one can construct a sequence of Galerkin approximations $u_{bl,n,m}$  by Schauder's fixed point theorem. Then, as the estimate  \eqref{aprioriestimate} holds for $u_{bl,n,m}$ uniformly in $m$ and $n$,  the sequence $D u_{bl,n,m}$ is bounded in $L^p(\Omega_{bl,n})$ uniformly in $m$. Sending $m$ to infinity yields a solution $u_{bl,n}$, the convergence of the nonlinear stress tensor follows from Minty's trick. Note that one can then  perform on $u_{bl,n}$ the manipulations of the previous paragraph, so that it satisfies \eqref{aprioriestimate} uniformly in $n$.

 \medskip
 {\em Step 2}: we let  $n$ go to infinity. We first extend $u_{bl,n}$ by $0$ for $y_2 > n$, and fix  $M > |A|$.  From the uniform estimate \eqref{aprioriestimate}, we get easily  the following convergences (up to a subsequence):  
 \begin{equation}
 \begin{aligned}
 & u_{bl,n} \rightarrow u_{bl}  \mbox{ weakly in } W^{1,p}_{loc}(\Omega_{bl}), \\
 & D u_{bl,n} \rightarrow D u_{bl}  \mbox{ weakly in } L^p(\Omega^-_{bl}), \\
 & D u_{bl,n} {\bbbone}_{|Du_{bl,n}| < M}\rightarrow V_1 \mbox{ weakly in } L^2(\Omega^+_{bl}), \quad \mbox{ weakly-* in }  L^\infty(\Omega_{bl}^+), \\
 & D u_{bl,n} {\bbbone}_{|D u_{bl,n}| \ge M}\rightarrow V_2  \mbox{ weakly in } L^p(\Omega^+_{bl}). 
\end{aligned} 
 \end{equation}
 Of course, $D u_{bl} = V_1 + V_2$ in $\Omega_{bl}^+$. A key point is that  
 $$\mbox{$S(A + D u_{bl,n}) - S(A)$ is bounded uniformly in $n$ in $(L^{p}(\Omega^+_{bl}))' =  L^{p'}(\Omega^+_{bl})$ and in $\left(L^2(\Omega^+_{bl}) \cap L^\infty(\Omega_{bl}^+)\right)'$.}$$  
 and converges weakly-* to some $S^+$ in that space. To establish this uniform bound, we treat separately   
 $$ S_{n,1} \: :=  \: (S(A + D u_{bl,n}) - S(A)) {\bbbone}_{|Du_{bl,n}| < M}, \quad  S_{n,2} \: :=  \: (S(A + D u_{bl,n}) - S(A)) {\bbbone}_{|Du_{bl,n}| \ge M}. $$
\begin{itemize}
\item  For $S_{n,1}$, we use  the inequality \eqref{ineq3}. It gives $|S_{n,1}| \le C |D u_{bl,n}| {\bbbone}_{|Du_{bl,n}| < M}$, 
 which provides a uniform bound in  $L^{2} \cap L^\infty$, and so in particular in $L^{p'}$ and in $L^2$. 
\item For $S_{n,2}$, we use first that $|S_{n,2}| \le C |D u_{bl,n}|^{p-1} {\bbbone}_{|D u_{bl,n}| \ge M}$, so that it is uniformly bounded in $L^{p'}$. We use then \eqref{ineq3}, so that 
$ |S_{n,2}| \le C |D u_{bl,n}| {\bbbone}_{|D u_{bl,n}| \ge M}$, which yields a uniform bound in $L^p$, in particular in $(L^2 \cap L^\infty)'$ ($p \in ]1,2]$).
\end{itemize}
From there, standard manipulations give 
$$ \int_{\Omega_{bl}^+} (S(A+ D u_{bl,n}) - S(A)) : D u_{bl,n} \rightarrow \int_{\Omega_{bl}^+} S^+ : (V_1 + V_2) = \int_{\Omega_{bl}^+} S^+ : D u_{bl} $$
One has even more directly 
$$ \int_{\Omega_{bl}^-}  S(D u_{bl,n}) : D u_{bl,n}  \rightarrow \int_{\Omega_{bl}^-} S^- : D u_{bl} $$
and one concludes by Minty's trick that $S^+ = S(A+D u_{bl}) - S(A)$, $\: S^- = S(D u_{bl})$. It follows that  $u_{bl}$ satisfies  \eqref{BL1}-\eqref{BL2} in a weak sense. Finally, one can  perform on $u_{bl}$ the manipulations of the previous paragraph, so that it satisfies \eqref{aprioriestimate}.

 \subsubsection*{Uniqueness}

 Let $u_{bl}^1$ and $u_{bl}^2$ be weak solutions of \eqref{BL1}-\eqref{BL2}, that is satisfying the variational formulation 
 \begin{equation} \label{VF}
 \int_{\Omega_{bl}^+} S(A+ D u_{bl}^i) : D \varphi \: + \: \int_{\Omega_{bl}^-} S(D u_{bl}^i) : D \varphi =  -\int_{y_2 = 0} S(A)n \cdot \varphi {\rm\,d}S, \quad i=1,2 
 \end{equation}
  for all smooth divergence free fields $\varphi \in C^\infty_c(\Omega_{bl})$.  
 The point is then to replace $\varphi$ by $u_{bl}^1 - u_{bl}^2$, to obtain  
 \begin{equation}  \label{zeroidentity}
 \int_{\Omega_{bl}^+} \left( S(A+ D u_{bl}^1) - S(A + D u_{bl}^2) \right) : D (u_{bl}^1 - u_{bl}^2) \: + \: \int_{\Omega_{bl}^-} (S(D u_{bl}^1) - S(D u_{bl}^2) :  D (u_{bl}^1 - u_{bl}^2) = 0 . \end{equation}
Rigorously, one constructs by convolution a sequence $\varphi_n$ such that $D \varphi_n$ converges appropriately to $D u_{bl}^1 - D u_{bl}^2$. In the case $p < 2$, the convergence holds in 
$(L^2(\Omega_{bl}^+) \cap L^\infty(\Omega_{bl}^+))  + L^p(\Omega_{bl}^+)$, respectively in $L^p(\Omega_{bl}^-)$. One can pass to the limit as $n$ goes to infinity because  
$$ S(A+D u_{bl}^1) - S(A+D u_{bl}^2) =    \left( S(A+D u_{bl}^1) - S(A)  \right) + \left( S(A) - S(A+D u_{bl}^2) \right), $$
respectively $S(D u_{bl}^1) - S(D u_{bl}^2)$, belongs to the dual space: see the arguments of the construction scheme of section \ref{sectionWP}. 
Eventually, by strict convexity of  $M \rightarrow |M|^p$ ($p > 1$), \eqref{zeroidentity} implies that  $D u_{bl}^1 = D u_{bl}^{2}$. This implies that $u_{bl}^1 - u_{bl}^2$ is a constant (dimension is $2$),  and due to the zero boundary condition at $\pa \Omega_{bl}$, we get $u_{bl}^1 = u_{bl}^2$. 

\subsection{Saint-Venant estimate}  \label{subsecstvenant}
We focus in this paragraph on the asymptotic behaviour of $u_{bl}$ as $y_2$ goes to infinity. The point is to show exponential convergence of $u_{bl}$ to a constant field.  At first, we can use interior regularity results for the generalized Stokes equation in two dimensions. In particular, pondering on the results of \cite{Wolf} for $p < 2$, and \cite{Kaplicky2002} for $p \ge 2$, we have : 
\begin{Lemma}  \label{lemma_unifbound}
The solution built in Theorem \ref{thWP} satisfies: $u_{bl}$ has $C^{1,\alpha}$ regularity over $\Omega_{bl}  \cap \{ y_2 > 1\}$ for some $0 < \alpha < 1$. In particular, $\na u_{bl}$ is bounded uniformly over $\Omega_{bl}  \cap \{ y_2 > 1\}$.
\end{Lemma}
 {\em Proof}. Let $0 \le  t < s$. We define $\Omega_{bl}^{t,s} \: := \:  \Omega_{bl} \cap \{t < y_2 \le s\}$. Note that $\Omega_{bl} \cap \{y_2 > 1 \} = \cup_{t \in \N_*} \Omega_{bl}^{t,t+1}$. Moreover, from the {\it a priori estimate} \eqref{aprioriestimate} or \eqref{es:Dup}, we deduce easily that  
 \begin{equation} \label{uniformLp}
  \| D u_{bl} \|_{L^p(\Omega_{bl}^{t,t+2})} \le C 
 \end{equation}
  for all $t > 0$,  for a constant $C$ that does not depend on $t$. We then introduce: 
  \begin{equation*}
    v_t  \: := \:  2 A \left(0,y_2 - (t+\frac{1}{2})\right) + u_{bl} \: - \frac{1}{2} \: \int_{\Omega_{bl}^{t-\frac{1}{2},t+\frac{3}{2}}} u_{bl} \dy, \quad  y \in \Omega_{bl}^{t-\frac{1}{2},t+\frac{3}{2}}, \quad \forall t \in \N_*.  
  \end{equation*}
  From \eqref{BL1}: 
  $$    -\dive(S(D v_t)) + \nabla p_{bl} = 0, \quad \dive v_t= 0 \quad \mbox{ in } \: \Omega_{bl}^{t-1/2,t+3/2}, \quad \forall t \in \N_*. $$
 Moreover, we get for some $C$ independent of $t$: 
 \begin{equation} \label{uniformLpv}
    \| v_t \|_{W^{1,p}(\Omega_{bl}^{t-1/2,t+3/2})} \le C \quad  \forall t \in \N_*.
 \end{equation}
Note that this $W^{1,p}$ control follows from \eqref{uniformLp}: indeed, one can apply  the Poincar\'e inequality for functions with zero mean, and then the Korn inequality. One can then ponder on the interior regularity results of articles  \cite{Wolf} and \cite{Kaplicky2002}, depending on the value of $p$: $v_{t}$ has  $C^{1,\alpha}$ regularity over $\Omega_{bl}^{t,t+1}$ for some $\alpha \in (0,1)$ (independent of $t$):  for some $C'$,  
$$\|  v_t \|_{C^{1,\alpha}(\Omega_{bl}^{t,t+1})} \le C', \quad \mbox{and in particular} \quad   \| \na v_t \|_{L^\infty(\Omega_{bl}^{t,t+1})}  \le C'  \quad  \forall t \in \N_*.$$
Going back to $u_{bl}$ concludes the proof of the lemma. 

\medskip
We are now ready to establish a  keypoint in the proof of Theorem \ref{thEC}, called a   
{\it Saint-Venant estimate}: namely, we show that the energy of the solution located above $y_2 = t$ decays exponentially with $t$. In our context, a good energy is 
$$ E(t) \: :=  \:   \int_{\{y_2 > t\}} | \na u_{bl} |^2 \dy $$
for $t > 1$. Indeed, from Lemma \ref{lemma_unifbound}, there exists $M$ such that $|D u_{bl}| \le M$ for all $y$ with $y_2 > 1$. In particular, in the case $p < 2$,  when localized above  $y_2 =1$, the energy functional that appears at the left hand-side of \eqref{aprioriestimate} only involves the $L^2$ norm of the symmetric gradient (or of the gradient by the homogeneous Korn inequality, {\it cf} the appendix).  Hence, $\na u_{bl} \in L^2(\Omega_{bl} \cap \{ y_2 > 1 \})$. The same holds for $p \ge 2$, thanks to \eqref{es:Dup}.
\begin{Proposition} \label{expdecay}
There exists $C,\delta > 0$, such that $E(t) \le C \exp(-\delta t)$.
\end{Proposition}
{\em Proof}. Let $t > 1$, $\Omega_{bl}^t := \Omega_{bl} \cap \{ y_2 > t \}$. Let $M$ such that $|D u_{bl}|$ is bounded by $M$ over $\Omega_{bl}^1$, which exists due to Lemma~\ref{lemma_unifbound}.  As explained just above, one has   
$ \int_{\Omega_{bl}^1 }  |D u_{bl}|^2 < +\infty$, and by Korn inequality $E(1)$ is finite. In particular, $E(t)$ goes to zero as $t \rightarrow +\infty$ and the point is to quantify the speed of convergence. By the use of inequality   \eqref{ineq1} (with $a = A$, $b = A + D u_{bl}$), we find
\begin{equation} \label{boundE}
\begin{aligned} 
E(t) \le C \int_{\Omega_{bl}^t} |D u_{bl}|^2 \dy & \le  
 C' \int_{\Omega_{bl}^t}
 	 \left( |A + D u_{bl}|^{p-2}(A + D u_{bl}) - | A |^{p-2} A   \right): D u_{bl}  \dy
	 \\ &
		 \le  C' \lim\limits_{n \to \infty } \int_{\Omega_{bl}} 
	 	 \left( |A + D u_{bl}|^{p-2}( A + D u_{bl}) - |A |^{p-2} A   \right): D u_{bl} \,  \chi_n (y_2) \dy
	\end{aligned}
	\end{equation}
for a smooth $\chi_n$ with values in $[0,1]$ such  that  $\chi_n = 1$ over $[t,t+n$], $\chi_n=0$ outside $[t-1, t+n+1]$, and $|\chi'_n| \le 2$. Then,  we integrate by parts the right-hand side, taking into account the first equation in \eqref{BL1}. We write 
\begin{align} 
\nonumber
& \int_{\Omega_{bl}} 
 \left( |A + D u_{bl}|^{p-2}( A + D u_{bl}) - |A |^{p-2} A   \right): D u_{bl} \,  \chi_n (y_2) \dy  \\
\nonumber
= & - \int_{\Omega_{bl}} \na p_{bl} \cdot  u_{bl}  \chi_n(y_2) \dy - \int_{\Omega_{bl}} \left( (S(A + D u_{bl}) - S(A)) \left( \begin{smallmatrix} 0 \\ \chi'_n \end{smallmatrix} \right) \right) \cdot u_{bl} \dy  \\
\label{I1I2}
= &\int_{\Omega_{bl}} \left(S(A) - S(A+D u_{bl})) \left( \begin{smallmatrix} 0 \\ \chi'_n \end{smallmatrix} \right) \right) \cdot u_{bl} \dy +  \int_{\Omega_{bl}} p_{bl} \chi'_n  u_{bl,2} \dy     \: := \: I_1 + I_2. 
\end{align} 
To estimate $I_1$ and $I_2$, we shall make use of simple vector inequalities. Namely: 
\begin{equation}  \label{ineq3}
\mbox{for all $p \in ]1,2]$,   for all vectors $a,b$,  $a\neq 0$, }, \quad | |b|^{p-2} b - |a|^{p-2} a | \: \le \: C_{p,a} \, |b-a| ,
\end{equation} 
whereas
\begin{equation} \label{ineq3duo}
\mbox{for all $p > 2$,   for all vectors $a,b$, $|b| \le M$}, \quad | |b|^{p-2} b - |a|^{p-2} a | \: \le \: C_{p,a,M} \, |b-a| .
\end{equation}
The latter is a simple application of the finite increment inequality. As regards the former, we distinguish between two cases: 
\begin{itemize}
\item If $|b-a| < \frac{|a|}{2}$, it follows from the finite increments inequality. 
\item If $|b-a| \ge \frac{|a|}{2}$, we simply write 
\begin{align*}
| |b|^{p-2} b - |a|^{p-2} a | \: & \le \: |b|^{p-1} + |a|^{p-1} \: \le \: (3^{p-1} + 2^{p-1}) |b-a|^{p-1}   \le \: (3^{p-1} + 2^{p-1}) (\frac{|a|}{2})^{p-2} |b-a| 
\end{align*}
using that $\left(\frac{2 |b-a|}{|a|}\right)^{p-1} \le \left(\frac{2 |b-a|}{|a|}\right)$ for $1 < p \le 2$. 
\end{itemize}

\medskip
We shall also make use of the following: 
\begin{Lemma} \label{lem_averages} For any height $t > 0$
\begin{description}
\item[i)] $\int_{\{y_2 =t\}} u_{bl,2} = 0$. 
\item[ii)] $\int_{\{y_2=t\}} (S(A+D u_{bl}) - S(A)) \left( \begin{smallmatrix} 0 \\ 1 \end{smallmatrix}  \right) \cdot \left( \begin{smallmatrix} 1 \\ 0 \end{smallmatrix}  \right) = 0$. 
\end{description}
\end{Lemma}
{\em Proof of the lemma}.  

\smallskip
\noindent
i) The integration of the divergence-free condition over $\Omega_{bl}^{0,t}$ leads to
\begin{align*}
0 = \int_{\Omega_{bl}^{0,t}} \dive u_{bl} & = \int_{\{y_2 = t \}} u_{bl,2} -  \int_{\{y_2 = 0^+ \}}u_{bl,2} = \int_{\{y_2 = t \}} u_{bl,2} -  \int_{\{y_2 = 0^- \}} u_{bl,2} \\
& = \int_{\{y_2 = t \}} u_{bl,2} - \int_{\Omega_{bl^-}} \dive u_{bl} + \int_{\Gamma_{bl}} u_{bl} \cdot n = \int_{\{y_2 = t \}} u_{bl,2} ,
\end{align*}
where the second and fourth inequalities come respectively from the no-jump condition of $u_{bl}$ at $y_2 = 0$ and the Dirichlet condition at $\Gamma_{bl}$. 

\smallskip
\noindent
ii) By integration of the first equation in \eqref{BL1} over $\Omega_{bl}^{0,t}$ we get: 
$$ \int_{y_2 = t} (S(A+D u_{bl}) - S(A) - p_{bl} Id) \left( \begin{smallmatrix} 0 \\ 1 \end{smallmatrix}  \right) = 
 \int_{y_2 = 0^+} (S(A+D u_{bl}) - S(A) - p_{bl} Id) \left( \begin{smallmatrix} 0 \\ 1 \end{smallmatrix}  \right). $$
 In particular, the quantity 
$$ I := \int_{y_2 = t} (S(A+D u_{bl}) - S(A) - p_{bl} Id) \left( \begin{smallmatrix} 0 \\ 1 \end{smallmatrix}  \right)  \cdot \left( \begin{smallmatrix} 1 \\ 0 \end{smallmatrix}  \right) =  
 \int_{y_2 = t} ( S(A+D u_{bl}) - S(A)) \left( \begin{smallmatrix} 0 \\ 1 \end{smallmatrix}  \right)  \cdot \left( \begin{smallmatrix} 1 \\ 0 \end{smallmatrix}  \right)  $$
is independent of the variable $t$. To show that it is zero, we apply inequality \eqref{ineq3} or \eqref{ineq3duo} with $a = A$ and $b = A + D u_{bl}$, so that 
$$I^2 \: \le \:  C \left( \int_{\{y_2=t\}} |D u_{bl}| \right)^2 \: \le \: C'  \int_{\{y_2=t\}} |D u_{bl}|^2 $$
($C$ is bounded by Lemma~\ref{lemma_unifbound}). As $D u_{bl}$ belongs to $L^2(\Omega^1_{bl})$, there exists a sequence $t_n$ such that
$ \int_{\{y_2=t_n\}} |D u_{bl}|^2 \rightarrow 0$ as $n \rightarrow +\infty$. It follows that $I = 0$. 
This concludes the proof of the Lemma. 

\medskip
We can now come back to the treatment of $I_1$ and $I_2$. 
\begin{itemize}
\item Treatment of $I_1$. 
\end{itemize}
We note that $\chi'_n$ is supported in $[t-1,t] \cup [t+n,t+n+1]$.  By Lemma \ref{lem_averages} we can write 
\begin{align}\label{I1a}
 I_1 & =  \int_{\Omega_{bl}^{t-1,t}} \left(S(A) - S(A+D u_{bl})) \left( \begin{smallmatrix} 0 \\ \chi'_n \end{smallmatrix} \right) \right) \cdot (u_{bl} - \overline{c}) \\
 & + \:  \int_{\Omega_{bl}^{t+n,t+n+1}} \left(S(A) - S(A+D u_{bl})) \left( \begin{smallmatrix} 0 \\ \chi'_n \end{smallmatrix} \right) \right) \cdot (u_{bl} - \overline{c}_n) \: := I_{1,1} + I_{1,2} ,
 \end{align}
where 
\begin{equation}\label{mean}
\overline{c} := \dashint_{\Omega_{bl}^{t-1,t}} u_{bl} = 
 \left( \dashint_{\Omega_{bl}^{t-1,t}} u_{bl,1} , 0 \right) \quad \mbox{ and } \quad  
\overline{c}_n := \dashint_{\Omega_{bl}^{t+n,t+n+1}} u_{bl} = 
 \left( \dashint_{\Omega_{bl}^{t+n,t+n+1}} u_{bl,1} , 0 \right).
 \end{equation}
Again, we apply inequality \eqref{ineq3} or \eqref{ineq3duo} to find 
$$ I_{1,1} \le C \int_{\Omega_{bl}^{t-1,t}} | D u_{bl} | \, |u_{bl} - \overline{c}| $$
and by the Poincar\'e inequality for functions with zero mean, we easily deduce that 
$$ I_{1,1} \le C' \int_{\Omega_{bl}^{t-1,t}} | \na u_{bl} |^2  = C' \left( E(t-1) - E(t) \right). $$  
An upper bound on $I_{1,2}$ can be derived in the same way: 
$$ I_{1,2} \le C' (E(t+n) - E(t+n+1)) ,$$  
where the right-hand side going to zero as $n \rightarrow +\infty$ since $E(t') \to 0$ as $t'\to \infty$. Eventually: 
\begin{equation} \label{estimI1} 
\limsup_{n \rightarrow +\infty} I_{1} \: \le \: C \left(E(t-1) - E(t) \right). 
\end{equation}
\begin{itemize}
\item Treatment of $I_2$. 
\end{itemize}
We can again use the decomposition 
\begin{equation}\label{I2a}
 I_2 = \int_{\Omega_{bl}^{t-1,t}}  p_{bl} \chi'_n  u_{bl,2}   \: + \: \int_{\Omega_{bl}^{t+n,t+n+1}}  p_{bl} \chi'_n  u_{bl,2} \: := \: I_{2,1} + I_{2,2}. 
 \end{equation}
From Lemma \ref{lem_averages} i), we infer that 
$$ \int_{\Omega_{bl}^{t-1,t}} \chi'_n(y_2)  u_{bl,2}(y) \, \dy = 0. $$
 By standard results, there exists $w \in H^1_0(\Omega_{bl}^{t-1,t})$ satisfying $ \dive w(y) =  \chi'_n(y_2)  u_{bl,2}(y), \quad y \in  \Omega_{bl}^{t-1,t}$, and the estimate 
 $$\| w \|_{H^1(\Omega_{bl}^{t-1,t})} \le C \| \chi'_n(y_2)  u_{bl,2}(y) \|_{L^2(\Omega_{bl}^{t-1,t})} \le C'  \| u_{bl,2} \|_{L^2(\Omega_{bl}^{t-1,t})}, $$
 for constants $C,C'$ that do not depend on $t$. As $w$ is zero at the boundary: 
 $$ I_{1,2} = \int_{\Omega_{bl}^{t-1,t}}  p_{bl} \dive w = - \int_{\Omega_{bl}^{t-1,t}} \na  p_{bl} \cdot w =  \int_{\Omega_{bl}^{t-1,t}} (S(A+D u_{bl}) - S(A)) \cdot \na w $$
where the last equality comes from \eqref{BL1}.  We find as before ({\it cf} \eqref{ineq3} or \eqref{ineq3duo}): 
\begin{align*}
|I_{1,2}| & \le C \int_{\Omega_{bl}^{t-1,t}} | D u_{bl} | |\na w|  \le C \| D u_{bl} \|_{L^2(\Omega_{bl}^{t-1,t})} \| \na w \|_{L^2(\Omega_{bl}^{t-1,t})} \\
& \le  C' \| D u_{bl} \|_{L^2(\Omega_{bl}^{t-1,t})} \,  \|u_{bl,2} \|_{L^2(\Omega_{bl}^{t-1,t})} \le C''  \| \na u_{bl} \|_{L^2(\Omega_{bl}^{t-1,t})}^2 
\end{align*}
where we have controlled the $L^2$ norm of $u_{bl,2}$ by the $L^2$ norm of its gradient (we recall that $u_{bl,2}$ has zero mean).  A similar treatment can be performed with $I_{2,2}$, so that $I_{2,1} \le C (E(t-1) - E(t))$, $\: I_{2,2} \le C ( E(t+n) - E(t+n+1))$ and 
 \begin{equation} \label{estimI2} 
\limsup_{n \rightarrow +\infty} I_{2} \: \le \: C \left(E(t-1) - E(t) \right). 
\end{equation}

\medskip
Finally, combining \eqref{boundE}, \eqref{I1I2}, \eqref{estimI1} and \eqref{estimI2}, we get 
$$ E(t) \le C (E(t-1) - E(t))  $$
for some $C > 0$. It is well-known that this kind of differential inequality implies the exponential decay of Proposition \ref{expdecay} (see the appendix). The proof of the Proposition is therefore complete. We have now all the ingredients to show Theorem \ref{thEC}. 

\medskip
{\em Proof of Theorem \ref{thEC}}. 
Thanks to the regularity Lemma \ref{lemma_unifbound}, we know that $\na u_{bl}$ is uniformly bounded over $\Omega_{bl}^1$, and belongs to $L^2(\Omega_{bl}^1)$.  Of course, this implies that 
$\na u_{bl}$ belongs to $L^q(\Omega_{bl}^1)$ for all $q \in [2,+\infty]$. More precisely, combining the $L^\infty$ bound with the $L^2$ exponential decay of Proposition \ref{expdecay}, we have that 
\begin{equation} \label{expdecayLq}
  \| \na u_{bl} \|_{L^q(\Omega_{bl}^t)} \le C \exp(-\delta t)  
 \end{equation}
 (for some $C$ and $\delta$ depending on $q$). This exponential decay extends straightforwardly to all $1 \le q < +\infty$.  Let us  now fix $q > 2$. To understand the behavior of $u$ itself, we write the Sobolev inequality: for all $y$ and $y' \in B(y,r)$, 
 \begin{equation} \label{Sobolevineq}
 |u(y') - u(y)| \: \le \: C r^{1-\frac{2}{q}}  \left( \int_{B(y,2r)} |\na u(z)|^q dz \right)^{1/q}. 
 \end{equation}
 We deduce from there that: for all $y_2 > 2$, for all $s \ge 0$,  
\begin{align*}
& |u_{bl}(y_1,y_2+s) - u_{bl}(y_1,y_2)| \\
& \le |u_{bl}(y_1,y_2+s) - u_{bl}(y_1,y_2+ \lfloor s \rfloor) | + \sum_{k=0}^{\lfloor s \rfloor-1} |u_{bl}(y_1,y_2+k+1) - u_{bl}(y_1,y_2+k)| \\
  & \le C \left( \| \na u_{bl} \|_{L^q(B((y_1,y_2+s)^t,1))}  + \sum_{k=0}^{\lfloor s \rfloor-1}  \| \na u_{bl} \|_{L^q(B((y_1,y_2+k)^t,1))} \right) \\
  & \le  C' \left( e^{-\delta(y_2+s)} + \sum_{k=0}^{\lfloor s \rfloor-1} e^{-\delta (y_2 + k)} \right) 
  \end{align*}
 where the last inequality comes from \eqref{Sobolevineq}. This implies that $u_{bl}$ satisfies the Cauchy criterion uniformly in $y_1$, and thus converges uniformly in $y_1$ to some $u^\infty = u^\infty(y_1)$ as $y_2 \rightarrow +\infty$.  To show that $u^\infty$ is a constant field, we rely again on \eqref{Sobolevineq}, which yields for all $|y_1 - y'_1| \le 1$:
 $$ |u_{bl}(y_1,y_2) - u_{bl}(y'_1, y_2) | \le C |y_1 - y'_1|^{1 - \frac{2}{q}} \| \na u_{bl} \|_{L^q(B((y_1,y_2)^t,1))} \le  C' e^{-\delta y_2}. $$
Sending $y_2$ to infinity gives: $u^\infty(y_1) = u^\infty(y'_1)$. Finally, the fact that $u^\infty$ is a horizontal vector field follows from Lemma \ref{lem_averages}, point i). This concludes the proof of the Theorem~\ref{thEC}. 

\medskip
Eventually, for later purposes, we state  
\begin{Corollary} {\bf (higher order exponential decay)} \label{higherorder}
\begin{itemize}
\item There exists $\alpha \in (0,1)$, such that for all $s \in [0,\alpha)$, for all $1 \leq q < +\infty$,  one can find $C$ and $\delta > 0$ with  
$$ \| u_{bl} - u^\infty \|_{W^{s+1,q}(\Omega_{bl}^t)} \le C \exp(-\delta t), \quad \forall t \ge 1. $$
\item There exists $\alpha \in (0,1)$, such that for all $s \in [0,\alpha)$, for all $1 \leq q < +\infty$,  one can find $C$ and $\delta > 0$ with  
$$ \| p_{bl} - p^t \|_{W^{s,q}(\Omega_{bl}^{t,t+1})} \: \le \: C \exp(-\delta t), \quad \forall t \ge 1,  \quad \mbox{for some constant $p^t$}.$$ 
\end{itemize}
\end{Corollary}
{\em Proof of the corollary}. It was established above that 
$$ |u(y_1,y_2+s) - u(y_1,y_2)| \le C'  \left( e^{-\delta'(y_2+s)} + \sum_{k=0}^{\lfloor s \rfloor -1} e^{-\delta' (y_2 + k)}\right) . $$
for some $C'$ and $\delta' > 0$. From there, after sending $s$ to infinity, it is easily deduced that 
$$ \| u_{bl} - u^\infty \|_{L^q(\Omega_{bl}^t)} \le C \exp(-\delta t) .$$
It then remains to control the $W^{s,q}$ norm of $\na u_{bl}$. This control comes from the $C^{0,\alpha}$ uniform bound on $\na u_{bl}$ over $\Omega_{bl}^1$, see Lemma \ref{lemma_unifbound}. By Sobolev imbedding, it follows that 
$$ \| \na u_{bl} \|_{W^{s,q}(\Omega_{bl}^{t,t+1})} \le C, \quad \forall s \in [0,\alpha), \forall 1\leq q < +\infty  $$ 
uniformly in $t$. Interpolating this bound with the bound 
$\| \na u_{bl} \|_{L^q(\Omega_{bl}^{t,t+1})} \le C' \exp(-\delta' t)$ previously seen, we get  
$$  \| \na u_{bl} \|_{W^{s,q}(\Omega_{bl}^{t,t+1})} \le C'' \exp(-\delta'' t), \quad \forall s \in [0,\alpha), \forall 1\leq q < +\infty  .$$ 
The first inequality of the Lemma follows. 

\medskip
The second inequality, on the pressure $p_{bl}$, is derived from the one on $u_{bl}$. This derivation is somehow standard, and we do not detail it for the sake of brevity.

\section{Error estimates, wall Laws}

\subsection{Approximation by the Poiseuille flow.}

We now go back to our primitive system \eqref{EQ1}. A standard estimate on $u^\eps$ leads to
$$ \int_{\Omega^\eps} |D u^\eps|^p \: \le \: \int_{\Omega^\eps} e_1 \cdot u^\eps. $$
The Korn inequality implies that 
$$   \int_{\Omega^\eps} |\na  u^\eps|^p \: \le \: C \int_{\Omega^\eps} |D u^\eps|^p $$
for a constant $C$ independent of $\eps$: indeed, one can extend $u^\eps$ by $0$ for $x_2 < \eps \gamma(x_1/\eps)$ and apply the inequality on the square $\T \times [-1,1]$, {\it cf} the appendix. 
Also, by the Poincar\'e inequality:
$$  | \int_{\Omega^\eps} e_1 \cdot u^\eps | \le C \| u^\eps \|_{L^p(\Omega^\eps)} \le C'  \| \na u^\eps \|_{L^p(\Omega^\eps)}. $$
We find that 
\begin{equation} \label{basic_estimate}
\| u^\eps \|_{W^{1,p}(\Omega^\eps)} \le C.
\end{equation} In particular, it provides strong convergence of $u^\eps$ in $L^p$ by the Rellich theorem (up to a subsequence). As can be easily guessed, the limit of $u^0$ in $\Omega$ is the generalized Poiseuille flow $u^0$. One can even obtain an error estimate by a direct energy estimate of the difference (extending $u^0$ and $p^0$ by zero in $R^\eps$).  We focus on the case $1 < p \le 2$, and comment briefly the easier case $p \ge 2$ afterwards. We write  $\ue = \uz + \we$ and $p^\ep = p^0 + q^\ep$. We find, taking into account \eqref{EQ2}: 
	\begin{equation}\label{EQ4}
	\begin{split}
	- \Div  \bS(\tD \uz + \tD \we) + \Div \bS(\tD \uz) + \nabla q^\ep & = {\bbbone}_{R^\eps}   e_1 \quad \mbox{ in } \Omega^\ep \setminus \Sigma_0 ,\\
	\Div \we & = 0 \quad \mbox{ in } \Omega^\ep , \\
	 \we & = 0 \quad \mbox{ on } \Gamma^\ep \cup \Sigma_1 ,  \\
	 \we & \mbox{ is periodic in } x_1 \mbox{ with period } 1 , \\
	[\we]|_{\Sigma_0} = 0,  \quad  [\tS ( \tD \uz + \tD \we ) {n}  - \tS ( \tD \uz ) {n}  - q^\ep n ]|_{\Sigma_0} & =  -\tS ( \tD \uz) {n}|_{x_2=0^+} .
	 \end{split}
	\end{equation}
In particular, performing an energy estimate and distinguishing between $\Omega$ and $R^\eps$, we find 
	\begin{equation}\label{weakEQ4}
	  \int_{\Omega} \left( \tS( \tD \uz  + \tD \we) - \tS( \tD \uz) \right) : \tD \we  + \int_{R^\ep} \tS( \tD \we ) : \tD \we  
	  =  - \int_{\Sigma_0} \tS (\tD \uz ){n}\vert_{x_2 = 0^+} \cdot \we {\rm d}S + \int_{R^\eps} e_1 \cdot \we  
	\end{equation}
Relying on  inequalities \eqref{ineq1}-\eqref{ineq2}, we get for any $M > \| Du^0 \|_{L^\infty}$:    
\begin{multline}\label{zeroap1}
	\| D \we \|_{L^p(\Omega \cap \{ |D \we| \ge M\})}^{p}   + \| D \we \|_{L^2(\Omega \cap \{ |D \we| \le M\})}^{2}  + \| D \we \|_{L^p(R^\eps)}^{p} \\ 
	\le C \left( \left| \int_{\Sigma_0}  \bS(\tD \uz)  {n}  \cdot \we  {\rm\,d}S \right| + \left| \int_{R^\eps} e_1 \cdot \we \right| \right)
	\end{multline}
Then by the H\"older inequality and  by Proposition~\ref{rescaledTracePoincare} in the appendix,  we have  that
	\begin{equation}\label{es1}
	| \int_{R^\ep} e_1 \cdot \we  | \leq  \eps^{\frac{p-1}{p}} \|   \we \|_{L^p(R^\eps)} 
	\leq C  \eps^{1 + \frac{p-1}{p}}  \| \nabla \we \|_{L^p(R^\ep)} .
	\end{equation}
Next, since $D \uz$ is given explicitly and uniformly bounded, the Proposition~\ref{rescaledTracePoincare} provides
	\begin{equation}\label{es2}
	\begin{split}
	& | \int_{\Sigma_0}  ( \tS(\tD \uz) {n}\vert_{x_2 = 0^+} \cdot   \we  {\rm\,d}S |
	 \leq C  \| \we \|_{L^{p}(\Sigma_0)} 
	 \leq C'   \ep^{\frac{p-1}{p}}  \| \nabla \we \|_{L^{p}(R^\ep)} .
	 \end{split}
	\end{equation}
Note that, as  $\we$ is zero at the lower boundary of the channel, we can extend it by $0$ below $R^\eps$ and apply Korn inequalities in a strip (see the appendix).  We find 
$$ \| \nabla \we \|_{L^p(R^\ep)}  \le C \| D \we \|_{L^p(R^\ep)}  $$
for some constant $C > 0$ independent of $\eps$. Summarising, we get

\begin{equation*}
	\| D \we \|_{L^p(\Omega \cap \{ |D \we| \ge M\})}^{p}   + \| D \we \|_{L^2(\Omega \cap \{ |D \we| \le M\})}^{2}  + \| D \we \|_{L^p(R^\eps)}^{p}  
	\le  C \ep^{\frac{p-1}{p}} \| D \we \|_{L^p(R^\ep)} 
\end{equation*}
and consequently  
\begin{equation}\label{zeroap2}
\| D \we \|_{L^p(\Omega \cap \{ |D \we| \ge M\})}^{p}   + \| D \we \|_{L^2(\Omega \cap \{ |D \we| \le M\})}^{2}  + \| D \we \|_{L^p(R^\eps)}^{p}  
	\le  C \ep
\end{equation}  
In the case $p \ge 2$ one needs to use \eqref{abp_1} instead of \eqref{ineq1}-\eqref{ineq2} what  yields 
\begin{equation} \label{zeroap3}
 \| D \we \|_{L^p(\Omega^\eps)} \le C \eps^{\frac{1}{p}}, \quad p \in [2, \infty). 
\end{equation}
\subsection{Construction of a refined approximation}
The aim of this section is to design a better approximation of the exact solution $u^\eps$ of \eqref{EQ1}. This approximation will of course involve the boundary layer profile $u_{bl}$ studied in the previous section. Consequences of this approximation in terms of wall laws will be discussed in paragraph \ref{parag_wall_laws}.

\medskip
From the previous paragraph, we know that the Poiseuille flow $u^0$ is the limit of $u^\eps$ in $W^{1,p}(\Omega)$. 
However, the extension of $u^0$ by $0$ in the rough part of the channel was responsible for a jump of the stress tensor at $\Sigma_0$. This jump was the main limitation of the error estimates \eqref{zeroap2}-\eqref{zeroap3},  and the reason for the introduction of the boundary layer term $u_{bl}$. Hence, we hope to have a better approximation replacing $u^0$ by $u^0 + \eps u_{bl}(\cdot/\eps)$. Actually, one can still improve the approximation, accounting for the so-called boundary layer tail $u^\infty$. More precisely, {\em in the Newtonian case}, a good idea is to replace $u^0$ by the solution  $u^{0,\eps}$ of the Couette problem:
$$ -\Delta u^{0,\eps} + \na p^{0,\eps} = 0, \quad \dive u^{0,\eps} = 0, \quad u^{0,\eps}\vert_{\Sigma_0} = \eps u^\infty, \quad  u^{0,\eps}\vert_{x_2 = 1} = 0.   $$
One then defines: 
$$ u^{\eps} = u^{0,\eps} + \eps (u_{bl}(\cdot/\eps) - u_\infty) + r^\eps \: \mbox{ in $\Omega$}, \quad u^\eps =   \eps u_{bl}(\cdot/\eps) \: \mbox{ in $R^\eps$},$$
where $r^\eps$ is a small divergence-free remainder correcting the $O(\exp(-\delta/\eps))$ trace of $u_{bl} - u^\infty$ at $\{x_2 = 1 \}$. 

\medskip
However, for technical reasons, the above approximation is not so successful in our context, so that we need to modify it a little. We proceed as follows. Let $N$ a large constant to be fixed later. We introduce: 
$$ \Omega^{\eps}_N := \Omega^\eps \cap \{x_2 > N \eps |\ln \eps| \}, \quad      \Omega^{\eps}_{0,N} = \Omega^\eps \cap \{ 0 < x_2 <  N \eps |\ln \eps| \}, \quad \mbox{and} \quad \Sigma_N = \Pi \times  \{ x_2 = N \eps  |\ln \eps| \}. $$ 
First, we introduce the solution $u^{0,\eps}$ of 
\begin{equation} \label{u0eps}
\left\{
  \begin{aligned}
  - \dive S(D u^{0,\eps}) + \na p^{0,\eps}  & = e_1, \quad x \in \Omega^\eps_N, \\
  \dive u^{0,\eps} & = 0, \quad x \in \Omega^\eps_N, \\
  u^{0,\eps}\vert_{\Sigma_N} & = \left( x \rightarrow  \left( \begin{smallmatrix} U'(0) x_2 \\ 0 \end{smallmatrix} \right) + \eps u^\infty \right)\vert_{\Sigma_N},  \\
    u^{0,\eps}\vert_{\{ x_2 = 1 \}} & = 0. 
  \end{aligned}
  \right.
  \end{equation}
  As for the generalized Poiseuille flow, the pressure $p^{0,\eps}$ is zero, and one has an explicit expression for $u^{0,\eps} = (U^\eps(x_2),0)$. In particular, one can check that 
\begin{equation} \label{explicit1} 
U^\eps(x_2) = \beta(\eps) -  \frac{(\sqrt{2})^{p'}}{p'} \left| \frac{1}{2} + \alpha(\eps) -x_2 \right|^{p'}, 
\end{equation}
where $\alpha(\eps)$ satisfies the equation ($x_{2,N} := N \eps |\ln \eps|$)
\begin{equation} \label{explicit2}
-\frac{1}{p'}(\sqrt{2})^{p'} \left( \left| \frac{1}{2} + \alpha(\eps) - x_{2,N} \right|^{p'} -  \left| \frac{1}{2} - \alpha(\eps)  \right|^{p'} \right) =  U'(0) x_{2,N}  + \eps U^\infty
\end{equation} 
and 
$$ \beta(\eps) = \frac{(\sqrt{2})^{p'}}{p'} \left| \frac{1}{2} - \alpha(\eps)\right|^{p'} . $$
 By the Taylor expansion, we find that 
  \begin{equation} \label{explicit3}
   \alpha(\eps) = - \sqrt{2}^{p'-4} \eps U^\infty + O(\eps^2 |\ln \eps|^2). 
  \end{equation}
  This will be used later.

  \medskip 
  Then, we consider the Bogovski problem 
  \begin{equation} \label{Bogov}
  \left\{
  \begin{aligned}
  \dive r^\eps & = 0 \quad \mbox{in} \:  \Omega^\eps_{N}, \\
   r^\eps\vert_{\Sigma_N} & = \eps (u_{bl}(\cdot/\eps) - u^\infty)\vert_{\Sigma_N}, \\
   r^\eps\vert_{\{x_2 = 1\}} & = 0. 
   \end{aligned}
   \right.
   \end{equation}
 Since $u^\infty = (U^\infty,0)$, note that 
   $$  \int_{\Sigma_N} \eps (u_{bl}(\cdot/\eps) - u^\infty) \cdot e_2 = \int_{\Omega^\eps_N \cup \overline{R^\eps}} {\rm div}_y u_{bl}(\cdot/\eps) = 0. $$
   Hence, the compatibility condition for solvability of \eqref{Bogov} is fulfilled: there exists  a solution 
   $r^\eps$ satisfying 
   $$ \| r^\eps \|_{W^{1,p}(\Omega^\eps_N)} \le C \eps \| u_{bl}(\cdot/\eps) - u^\infty \|_{W^{1-\frac{1}{p},p}(\Sigma_N)}. $$
Using the first estimate of Corollary \ref{higherorder}, we find 
\begin{equation} \label{estimreps}   
 \| r^\eps \|_{W^{1,p}(\Omega^\eps_N)} \le C \eps^{\frac{1}{p}} \exp(-\delta N |\ln \eps|). 
 \end{equation}
 
 \medskip
 Finally, we define the approximation $(u^\eps_{app}, p^\eps_{app})$ by the formula
 \begin{equation}   \label{uepsapp}
 u^\eps_{app}(x) = \left\{ 
\begin{aligned}  
& u^{0,\eps}(x) + r^\eps(x) \quad x \in \Omega^\eps_N ,\\
&  \left( \begin{smallmatrix} U'(0)x_2 \\ 0 \end{smallmatrix} \right) + \eps u_{bl}(x/\eps), \quad x \in  \Omega^\eps_{0,N}, \\
& \eps u_{bl}(x/\eps), \quad x \in  R^\eps ,
\end{aligned}
\right.
\end{equation}
whereas 
 \begin{equation}  
 p^\eps_{app}(x) = \left\{ 
\begin{aligned}  
& 0 \quad x \in \Omega^\eps_N ,\\
& p_{bl}(x/\eps) \quad \quad x \in  \Omega^\eps_{0,N} \cup R^\eps.
\end{aligned}
\right.
\end{equation}
With such a choice:   
$$u^\eps_{app}\vert_{\partial \Omega^\eps} = 0, \quad   \dive u^\eps_{app} = 0 \quad \mbox{over $\Omega^\eps_N \cup  \Omega^\eps_{0,N} \cup R^\eps$}.$$
Moreover, $u^\eps_{app}$ has zero jump at the interfaces $\Sigma_0$ and $\Sigma_N$: 
$$ [u^\eps_{app}]\vert_{\Sigma_0} = 0, \quad [u^\eps_{app}]\vert_{\Sigma_N} = 0. $$
Still, the stress tensor has a jump. More precisely, we find 
\begin{equation} \label{stressjump} 
\begin{aligned}
\left[S(D u^\eps_{app})n - p^\eps_{app}n\right]\vert_{\Sigma_0} & = 0, \\
\left[S(D u^\eps_{app})n - p^\eps_{app}n\right]\vert_{\Sigma_N} & = \left( S(D u^0_\eps + D r^\eps)\vert_{\{x_2 = (N \eps |\ln \eps|)^+\}} - S(A + Du_{bl}(\cdot/\eps))\vert_{\{x_2 = (N \eps |\ln \eps|)^-\}}\right) e_2 \\
& - p_{bl}(\cdot/\eps)  \vert_{\{x_2 = (N \eps |\ln \eps|)^-\}} e_2 .
\end{aligned}
\end{equation}
The next step is to obtain error estimates on $u^\eps - u^\eps_{app}$. 
\subsection{Error estimates}
We prove here: 
\begin{Theorem} {\bf (Error estimates)} \label{thmerror}
\begin{itemize}
\item For $1 < p \le 2$, there exists $C$ such that 
$$ \| u^\eps - u^\eps_{app} \|_{W^{1,p}(\Omega^\eps)} \le C (\eps |\ln \eps|)^{1+\frac{1}{p'}} .$$
\item For $p \ge 2$, there exists $C$ such that 
$$ \| u^\eps - u^\eps_{app} \|_{W^{1,p}(\Omega^\eps)} \le C (\eps |\ln \eps|)^{\frac{1}{p-1}+\frac{1}{p}} . $$
\end{itemize}
\end{Theorem}
\begin{Remark}
A more careful treatment would allow to get rid of the $\ln$ factor in the last estimate ($p \ge 2$). We do not detail this point here, as we prefer  to provide a unified treatment. Also, we recall that  the shear thinning case  ($1 < p \le 2$) has a much broader range of applications. More comments will be made on the estimates in the last paragraph \ref{parag_wall_laws}.   
\end{Remark}
 
{\em Proof of the theorem.} We write $v^\eps = u^\eps - u^\eps_{app}$, $q^\eps = p^\eps - p^\eps_{app}$. We start from the equation  
\begin{equation} \label{eq_error}
-\dive S(D u^\eps) + \dive S(D u^\eps_{app})  + \na q^\eps = e_1 + \dive S(D u^\eps_{app}) + \na p^\eps_{app} := F^\eps 
\end{equation}
satisfied in $\Omega^\eps \setminus (\Sigma_0 \cup \Sigma_N)$. 
 A quick computation shows that 
\begin{equation*}
F^\eps = 
\left\{
\begin{aligned}
& \dive S(D u^{0,\eps} + D r^\eps) - S(D u^{0,\eps}), \quad x \in \Omega^\eps_N, \\
& e_1, \quad x \in \Omega^\eps_{0,N} \cup R^\eps. 
\end{aligned}
\right.
\end{equation*} 
Defining 
$$  \langle F^\eps, v^\eps \rangle := \int_{\Omega^\eps_N} F^\eps \cdot v^\eps +  \int_{\Omega^\eps_{0,N}} F^\eps \cdot v^\eps + \int_{R^\eps} F^\eps \cdot v^\eps $$
we get:
$$ | \langle F^\eps, v^\eps \rangle | \le \alpha_\eps \| \na v^\eps \|_{L^p(\Omega^\eps_N)} + \beta_\eps \| v^\eps \|_{L^p(\Sigma_N)} + \| v^\eps \|_{L^1(\Omega^\eps \setminus \Omega^\eps_N)} $$
where 
$$ \alpha_\eps := \| S(D u^{0,\eps} + D r^\eps) - S(D u^{0,\eps}) \|_{L^{p'}(\Omega^\eps_N)}, \quad \beta_\eps := \| \left(S(D u^{0,\eps} + D r^\eps) - S(D u^{0,\eps})\right)e_2 \|_{L^{p'}(\Sigma_N)}. $$
We then use the inequalities 
\begin{equation} \label{poincarelike}
\begin{aligned}
  \| v^\eps \|_{L^p(\Sigma_N)} \: \le \: C (\eps |\ln \eps|)^{1/p'}  \| \na v^\eps \|_{L^p(\Omega^\eps)}, \\
  \| v^\eps \|_{L^1(\Omega^\eps \setminus \Omega^\eps_N)} \le C \eps^{\frac{1}{p'}}  \| v^\eps \|_{L^p(\Omega^\eps \setminus \Omega^\eps_N)}  \: \le \: C \eps^{\frac{1}{p'}} (\eps |\ln \eps|) \| \na v^\eps \|_{L^p(\Omega^\eps)}
  \end{aligned}
  \end{equation}
  (see the appendix for similar ones). We end up with 
\begin{equation}
   | \langle F^\eps, v^\eps \rangle | \le C \left( \alpha_\eps + \beta_\eps  (\eps |\ln \eps|)^{1/p'}  +  \eps^{\frac{1}{p'}} (\eps |\ln \eps|) \right) \| \na v^\eps \|_{L^p(\Omega^\eps)} .
\end{equation}

\medskip
Back to \eqref{eq_error}, after multiplication by $v^\eps$ and integration over $\Omega^\eps$, we find:
\begin{equation*}
\begin{aligned}
 & \int_{\Omega^\eps} \left( S(D u^\eps) - S(Du^\eps_{app}) \right) :\na v^\eps \\
& \le C \left( \alpha_\eps + \beta_\eps  (\eps |\ln \eps|)^{1/p'}  +  \eps^{\frac{1}{p'}} (\eps |\ln \eps|)  \right) \| \na v^\eps \|_{L^p(\Omega^\eps)} + \int_{\Sigma_N} \left( [S(D u^\eps_{app}) e_2]\vert_{\Sigma_N} \cdot v^\eps - [p^\eps_{app}]\vert_{\Sigma_N} v^\eps_2 \right). 
\end{aligned}
\end{equation*}
Let $p^{\eps,N}$ be a constant to be fixed later. As $v^\eps$ is divergence-free and zero at $\Gamma^\eps$, its flux through $\Sigma_N$ is zero:  $\int_{\Sigma_N} v^\eps_2 = 0$. Hence, we can add $p^{\eps,N}$ to the pressure jump $[p^\eps_{app}]\vert_{\Sigma_N}$ without changing the surface integral. We get:   

\begin{equation} \label{final_estimate}
\begin{aligned}
 & \int_{\Omega^\eps} \left( S(D u^\eps) - S(Du^\eps_{app}) \right) : \na v^\eps \\
& \le C \left( \alpha_\eps + \beta_\eps  (\eps |\ln \eps|)^{1/p'}  +  \eps^{\frac{1}{p'}} (\eps |\ln \eps|) \right) \| \na v^\eps \|_{L^p(\Omega^\eps)} + \int_{\Sigma_N} \left( [S(D u^\eps_{app}) e_2]\vert_{\Sigma_N} \cdot v^\eps - ([p^\eps_{app}]\vert_{\Sigma_N} - p^{\eps,N}) v^\eps_2 \right) \\
& \le \left( \alpha_\eps + \beta_\eps  (\eps |\ln \eps|)^{1/p'}  +  \eps^{\frac{1}{p'}} (\eps |\ln \eps|)  \right) \| \na v^\eps \|_{L^p(\Omega^\eps)} + \gamma_\eps \| v^\eps \|_{L^p(\Sigma_N)} \\
&  \le C  \left( \alpha_\eps + (\beta_\eps + \gamma^\eps)  (\eps |\ln \eps|)^{1/p'}  +  \eps^{\frac{1}{p'}} (\eps |\ln \eps|)  \right)  \| \na v^\eps \|_{L^p(\Omega^\eps)} ,
\end{aligned}
\end{equation}
where 
$$ \gamma_\eps := \| [\left(S(D u^\eps_{app})]\vert_{\Sigma_N} - ([p^\eps_{app}]\vert_{\Sigma_N} - p^{\eps,N}\right) e_2) \|_{L^{p'}(\Sigma_N)}.  
$$ 
Note that  we used again the first bound in \eqref{poincarelike} to go from the third to the fourth inequality. 
\begin{Lemma} \label{bounds}
For $N$ large enough, and a good choice of $p^{\eps,N}$ there exists $C = C(N)$ such that
$$ \alpha_\eps \le C \eps^{10}, \quad \beta^\eps \le C \eps^{10}, \quad \gamma_\eps \le C \eps  |\ln \eps|. $$
\end{Lemma}

\medskip
Let us temporarily admit this lemma. Then, we can conclude the proof of the error estimates:
\begin{itemize}
\item In the case $1 \le p \le 2$, we rely on the inequality established in \cite[Proposition~5.2]{GM1975}: for all $p \in ]1,2]$, there exists $c$ such that  for all $u,u' \in W_0^{1,p}(\Omega^\eps)$ 
$$  \int_{\Omega^\eps} \left( S(D u) - S(D u') \right) \cdot \na (u - u') \ge c \frac{\| Du - Du' \|^2_{L^p(\Omega^\eps)}}{(\| D u \|_{L^p(\Omega^\eps)} + \| D u' \|_{L^p(\Omega^\eps)})^{2-p}}  $$
We use this inequality with $u = u^\eps$, $u' = u^\eps_{app}$. With the estimate \eqref{basic_estimate} and the Korn inequality in mind, we obtain 
$$  \int_{\Omega^\eps} \left( S(D u^\eps) - S(Du^\eps_{app}) \right) \cdot \na v^\eps  \ge c \| \na v^\eps \|_{L^p}^2. $$
Combining this lower bound with the upper bounds on $\alpha_\eps, \beta_\eps, \gamma_\eps$ given by the lemma, we deduce from \eqref{final_estimate} the first error estimate in Theorem \ref{thmerror}. 
\item In the case $2 \le p$, we use the easier  inequality 
$$   \int_{\Omega^\eps} \left( S(D u) - S(D u') \right) \cdot \na (u - u') \ge c \| D u - D u' \|_{L^p(\Omega^\eps)}^p, $$ so that 
$$    \int_{\Omega^\eps} \left( S(D u^\eps) - S(Du^\eps_{app}) \right) \cdot \na v^\eps  \ge c \| \na v^\eps \|_{L^p(\Omega^\eps)}^p. $$
The second error estimate from Theorem \ref{thmerror} follows. 
\end{itemize}

\medskip 
The final step is to establish the bounds of Lemma \ref{bounds}. 

\medskip
{\em Bound on $\alpha_\eps$ and $\beta_\eps$}. From Corollary \ref{higherorder} and the trace theorem, we deduce that 
\begin{equation} \label{traceubl} 
\| u_{bl}(\cdot/\eps) - u^\infty \|_{W^{1+s-\frac{1}{q},q}(\{ x_2 = t \})} \le C \eps^{\frac{1}{q}-s-1}\exp(-\delta t/\eps) 
\end{equation}
for some $s < \alpha$ (where $\alpha \in (0,1)$) and any $q > \frac{1}{s}$. Let $q > \max(p',  \frac{2}{s})$. The solution $r^\eps$ of \eqref{Bogov} satisfies: 
$r^\eps \in W^{1+s,q}(\Omega^\eps_N)$ with 
\begin{equation*}
 \| r^\eps \|_{W^{1+s,q}(\Omega^\eps_N)} \le C  \eps^{\frac{1}{q}-s} \exp(-N\delta|\ln \eps|) 
 \end{equation*}
 so that by Sobolev imbedding
\begin{equation} \label{estimreps2}
\| D r^\eps \|_{L^\infty(\Sigma_N)} + \| D r^\eps \|_{L^q(\Sigma_N)}  + \| D r^\eps \|_{L^\infty(\Omega^\eps_N)}   \le  C \| D r^\eps \|_{W^{s,q}(\Omega^\eps_N)} \le C  \eps^{\frac{1}{q}-s} \exp(-N\delta|\ln \eps|) 
\end{equation}
This last inequality allows to evaluate $\beta_\eps$. Indeed, for $x \in \Sigma_N$, $C \ge |Du^{0,\eps}(x)| \ge c > 0$ uniformly in $x$. We can then use the upper bound \eqref{ineq3} for $p < 2$,  or \eqref{ineq3duo} for $p \ge 2$, to obtain 
\begin{equation}\label{betaep}
 \beta_\eps \le C \| D r^\eps \|_{L^{p'}(\Sigma_N)} \le C \| D r^\eps \|_{L^q(\Sigma_N)}  \le C'  \eps^{\frac{1}{q}-s} \exp(-N\delta|\ln \eps|) \le C' \eps^{10}\end{equation}
for $N$ large enough.

\medskip
To treat $\alpha_\eps$, we still have to pay attention to the cancellation of $D u^{0,\eps}$. Indeed, from the explicit expression of $u^{0,\eps}$, we know that there is some $x_2(\eps) \sim \frac{1}{2}$ at which $D u^{0,\eps}\vert_{x_2 = x_2(\eps)} = 0$.  Namely, we write 
\begin{align*}
& \int_{\Omega^\eps_N} |S(D u^{0,\eps} + D r^\eps) - S(D u^{0,\eps})|^{p'} \\
& = \int_{\{x \in \Omega^\eps_N,\, | x_2 - x_2(\eps) | \le \eps^{10 p'}\}}  |S(D u^{0,\eps} + D r^\eps) - S(D u^{0,\eps})|^{p'}    
 +   \int_{\{x\in \Omega^\eps_N,\, | x_2 - x_2(\eps) | \ge \eps^{10 p'}\}}  |S(D u^{0,\eps} + D r^\eps) - S(D u^{0,\eps})|^{p'}  \\
& := I_1 + I_2. 
\end{align*}
The first integral is bounded by
$$ I_1 \le C \int_{\{x\in \Omega^\eps_N,\, | x_2 - x_2(\eps) | \le \eps^{10 p'}\}} | D u^{0,\eps} |^p + | D r^\eps |^p \le C \eps^{10 p'}, $$
where we have used the uniform bound satisfied by $D u^{0,\eps}$ and $D r^\eps$ over $\Omega^\eps_N$, see \eqref{estimreps2}. 
For the second integral, we can distinguish between $p < 2$ and $p \ge 2$.  For $p < 2$, see \eqref{ineq3} and its proof, we get   
\begin{equation}\label{I_2}
 I_2 \le C  \int_{\{x\in \Omega^\eps_N,\, | x_2 - x_2(\eps) | \ge \eps^{10 p'}\}}  |D u^{0,\eps}|^{(p-2)p'} |D r^\eps|^{p'} 
\le C' \eps^{-M} \exp(-\delta'  N|\ln \eps|) 
\end{equation}
for some  $M, C', \delta' > 0$, see \eqref{estimreps2}. In the case $p \ge 2$, as  $D u^{0,\eps}$ and $D r^\eps$ are uniformly bounded, we derive a similar inequality by \eqref{ineq3duo}.  In both cases,  taking $N$ large enough, we obtain $I_2 \le C'' \eps^{10 p'}$, to end up with $\alpha_\eps \le C \eps^{10}$.  

\medskip
{\em Bound on $\gamma_\eps$}. We have 
\begin{align*}
 \gamma_\eps & 
  \le \| \left(S(D u^{0,\eps} + D r^\eps) - S(D u^{0,\eps})\right) e_2 \|_{L^{p'}(\Sigma_N)} 
 + \| \left(S(D u^{0,\eps}) - S(A)\right) e_2 \|_{L^{p'}(\Sigma_N)}  
 \\
 & +  \| (S(A)  - S(A + D u_{bl}(\cdot/\eps))) e_2 \|_{L^{p'}(\Sigma_N)}  + \| p_{bl}(\cdot/\eps) - p^{\eps,N} \|_{L^{p'}(\Sigma_N)} .
\end{align*}
The first term is $\beta_\eps$, so $O(\eps^{10})$  by previous calculations. The third term can be treated similarly to $\beta_\eps$. As $A \neq 0$, \eqref{ineq3} implies that 
\begin{equation}\label{I_3bis}
 \| (S(A)  - S(A + D u_{bl}(\cdot/\eps))) e_2 \|_{L^{p'}(\Sigma_N)}  \le C \| D u_{bl}(\cdot/\eps) \|_{L^{p'}(\Sigma_N)} \le C' \exp(-\delta' N \ln \eps) ,
 \end{equation}
where the last inequality can be deduced from \eqref{traceubl}. It is again $O(\eps^{10})$ for $N$ large enough.  For the second term of the right-hand side, we rely on the explicit expression of $u^{0,\eps}$.  On the basis of \eqref{explicit1}-\eqref{explicit3}, we find that 
$$ D(u^{0,\eps})\vert_{\Sigma_N} = A + O(\eps |\ln \eps|) $$ 
resulting in 
$$ \| \left(S(D u^{0,\eps}) - S(A)\right) e_2 \|_{L^{p'}(\Sigma_N)}  \le C \eps |\ln \eps|. $$
Finally, to handle the pressure term, we use the second term of Corollary \ref{higherorder}, which implies  
$$ \| p_{bl} - p^t \|_{L^q(\{ y_2 = t\})} \: \le \: C \exp(-\delta t)  \quad \mbox{for some constant $p^t$}.$$ 
We take $t  =  N \ln \eps$ and $p^{\eps,N} = p^t$ to get 
$$  \| p_{bl}(\cdot/\eps) - p^{\eps,N} \|_{L^{p'}(\Sigma_N)}  \le C' \exp(-\delta' N |\ln \eps|). $$
Taking $N$ large enough, we can make this term neglectible, say $O(\eps^{10})$. Gathering all contributions, we obtain $  \gamma_\eps \le C \eps |\ln \eps|$ as stated. 


\subsection{Comment on possible wall laws} \label{parag_wall_laws}
On the basis of the previous error estimates, we can now discuss the appropriate wall laws for a non-Newtonian flow above a rough wall. We focus here again on the shear thinning case  ($1 < p \le 2$).

\medskip
We first notice that the field $u^\eps_{app}$ (see  \eqref{uepsapp}) involves in a crucial way the solution $u^ {0,\eps}$ of \eqref{u0eps}. Indeed, we know from \eqref{estimreps} that  the contribution of $r^\eps$  in $W^{1,p}(\Omega^\eps_N)$ is very small  for $N$ large enough. Hence, the error estimate of Theorem \ref{thmerror} implies that   
$$ \| u^\eps - u^{0,\eps} \|_{W^{1,p}(\Omega^\eps_N)} = O((\eps |\ln \eps|)^{1+\frac{1}{p'}}) . $$
In other words, away from the boundary layer, $u^\eps$ is well approximated by $u^{0,\eps}$, with a power of $\eps$ strictly bigger than $1$. Although such estimate is unlikely to be optimal, it is enough to emphasize the role of the boundary layer tail $u^\infty$. Namely, the addition of the term  $\eps u^\infty$ in the Dirichlet condition for $u^{0,\eps}$ (see the third line of \eqref{u0eps}) allows to go beyond a $O(\eps)$ error estimate.  
{\it A contrario}, the  generalized Poiseuille flow $u^0$ leads to a $O(\eps)$ error only (away from the boundary layer). Notably, 
\begin{equation} \label{estimu0}
\| u^\eps - u^0 \|_{W^{1,p}(\Omega^\eps_N)} \ge \| u^{0,\eps}  - u^0 \|_{W^{1,p}(\Omega^\eps_N)}  - \| u^\eps -  u^{0,\eps} \|_{W^{1,p}(\Omega^\eps_N)} \ge c  \eps - o(\eps) \ge c' \eps ,
\end{equation}
where the lower bound for $u^{0,\eps}  - u^0$ is obtained using the explicit expressions.

\medskip
Let us further notice that instead of considering $u^{0,\eps}$, we could consider the solution  
$u^0_\eps$ of 
\begin{equation} \label{u0epsbis}
\left\{
\begin{aligned}
- \dive S(u^0_\eps)) + \na p^0_\eps  & = e_1, \quad x \in \Omega^\eps_N, \\
\dive u^0_\eps & = 0, \quad x \in \Omega^\eps_N, \\
u^0_\eps \vert_{\Sigma_0} & = \eps u^\infty,  \\
u^0_\eps\vert_{\{ x_2 = 1 \}} & = 0. 
\end{aligned}
\right.
\end{equation}
It reads $u^0_\eps = (U_\eps, 0)$ with 
$$ U_\eps(x_2) = \beta'(\eps) -  \frac{(\sqrt{2})^{p'}}{p'} \left| \frac{1}{2} + \alpha'(\eps) -x_2 \right|^{p'} $$
for $\alpha'$  and $\beta'$ satisfying 
$$ 
-\frac{1}{p'}(\sqrt{2})^{p'} \left( \left| \frac{1}{2} + \alpha'(\eps) \right|^{p'} -  \left| \frac{1}{2} - \alpha'(\eps)  \right|^{p'} \right) =   \eps u^{\infty}_1
\quad
\mbox{ and }
\quad 
 \beta'(\eps) = \frac{(\sqrt{2})^{p'}}{p'} \left| \frac{1}{2} - \alpha'(\eps)\right|^{p'}. $$
We can compare directly these expressions to \eqref{explicit1}-\eqref{explicit2} and deduce that 
$$ \| u^{0,\eps} - u^0_\eps \|_{W^{1,p}(\Omega^\eps_N)} = O(\eps |\ln \eps|),  $$
which in turn implies that 
\begin{equation} \label{estimu0eps}
\| u^\eps - u^0_\eps \|_{W^{1,p}(\Omega^\eps_N)} = O(\eps |\ln \eps|).  
\end{equation}
Hence, in view of \eqref{estimu0} and \eqref{estimu0eps}, we distinguish between two approximations (outside the boundary layer): 
\begin{itemize}
\item A crude approximation, involving the generalized Poiseuille flow  $u^0$. 
\item A refined approximation, involving $u^0_\eps$.
\end{itemize}
The first choice corresponds to the Dirichlet wall law $u\vert_{\Sigma_0} = 0$, and neglects the role of the roughness. The second choice takes it into account through the inhomogeneous Dirichlet condition: 
$u\vert_{\Sigma_0} = \eps u^\infty = \eps (U^\infty,0)$. Note that this last boundary condition can be expressed as  a wall law, although slightly abstract. Indeed, $U^\infty$ can be seen as a function of the tangential shear $(D(u^0)n)_\tau\vert_{\Sigma_0} = \pa_2 u^0_1\vert_{\Sigma_0} = U'(0)$, through the mapping    
$$ U'(0) \: \rightarrow \:  A := \left( \begin{smallmatrix} 0 & U'(0) \\ U'(0) & 0 \end{smallmatrix} \right) \:  \rightarrow \: u_{bl} \:\:  \mbox{solution of \eqref{BL1}-\eqref{BL2}}  \: \rightarrow \: U^\infty = \lim_{y_2 \rightarrow +\infty} u_{bl,1}. $$
Denoting by ${\cal F}$ this application, we write  
$$(u^0_\eps)_\tau \vert_{\Sigma_0} = \eps {\cal F}((D(u^0)n)_\tau\vert_{\Sigma_0}) \approx    \eps {\cal F}((D(u^0_\eps)n)_\tau\vert_{\Sigma_0})$$
whereas $ (u^0_\eps)_n  = 0$. This provides the following refined wall law : 
$$ u_n\vert_{\Sigma_0} = 0, \quad u_\tau\vert_{\Sigma_0} = \eps {\cal F}\bigl((D(u) n)_\tau\vert_{\Sigma_0}\bigr).  $$
 This wall law generalizes the Navier wall law derived in  the Newtonian case, where ${\cal F}$ is simply linear. Of course, it is not very explicit as it involves the nonlinear system \eqref{BL1}-\eqref{BL2}. 
 More studies will be necessary to obtain qualitative properties of the function ${\cal F}$, leading to a more  effective boundary condition.

\bigskip
{\bf Acknowledgements: } 
The work of AWK is partially  supported by Grant of National Science Center Sonata, No 2013/09/D/ST1/03692.

\section{Appendix : A few functional inequalities}
\begin{Proposition} {\bf (Korn inequality)}
Let $S_a  := \T \times (a, a+1)$, $a \in \R$. For all $1 < p < +\infty$, there exists $C > 0$ such that: for all $a \in \R$, for  all $u \in W^{1,p}(S_a)$, 
\begin{equation} \label{Korn1}
 \| \na u  \|_{L^p(S_a)} \: \le \: C \| D u \|_{L^p(S_a)}. 
\end{equation}
\end{Proposition}
{\em Proof.} Without loss of generality, we can show the inequality for $a = 0$: the independence of the constant $C$ with respect to $a$ follows from invariance by translation. Let us point out that the keypoint of the proposition is that the inequality is homogeneous. Indeed, it is well-known that the inhomogeneous Korn inequality 
\begin{equation}  \label{Korn2}
 \| \na u  \|_{L^p(S_0)} \: \le \: C' \left( \| D u \|_{L^p(S_0)} + \| u \|_{L^p(S_0)}\right)
\end{equation}
holds.  To prove the homogeneous one, we use reductio at absurdum : if \eqref{Korn1} is wrong, there exists a sequence $u_n$ in $W^{1,p}(S_0)$ such that 
\begin{equation} \label{absurd}
 \| \na u_n  \|_{L^p(S_0)} \: \ge \: n \| D u_n \|_{L^p(S_0)}. 
 \end{equation}
Up to replace $u_n$ by $u'_n := (u_n - \int_{S_0} u_n)/\| u_n\|_{L^p}$,  we can further assume that 
$$  \| u_n \|_{L^p} = 1, \quad \int_{S_0} u_n = 0.$$  
Combining \eqref{Korn2} and \eqref{absurd}, we deduce that $1 \ge \frac{n - C'}{C'} \| D(u_n) \|_{L^p}$ which shows that $D(u_n)$ converges to zero in $L^p$. Using again \eqref{Korn2}, we infer that $(u_n)$ is bounded in $W^{1,p}$, so that up to a subsequence it converges weakly  to some $u \in W^{1,p}$, with strong convergence in $L^p$ by Rellich Theorem. We have in particular  
\begin{equation} \label{contradiction}
\| u \|_{L^p} = \lim_n  \| u_n \|_{L^p} = 1, \quad \int_{S_0} u = \lim_n \int_{S_0} u_n = 0. 
\end{equation}
Moreover, as $D(u_n)$ goes to zero, we get $D(u) = 0$. This implies that $u$ must be a constant (dimension is 2), which makes the two statements of \eqref{contradiction} contradictory. 

\begin{Corollary}
Let $H_a := \T \times (a, + \infty)$. 
For all $1 < p < +\infty$, there exists $C > 0$ such that: for all $a \in \R$, for  all $u \in W^{1,p}(H_a)$, 
\begin{equation*} 
 \| \na u  \|_{L^p(H_a)} \: \le \: C \| D u \|_{L^p(H_a)}. 
\end{equation*}
\end{Corollary}
{\em Proof.}  From the previous inequality, we get for all $n \in \N$: 
$$ \int_{S_{a+n}} | \na u |^p \: \le C \:  \int_{S_{a+n}} | D u |^p. $$
The result follows by summing over $n$. 
\begin{Corollary}
Let $1 < p < +\infty$. There exists $C > 0$, such that for all $u \in W^{1,p}(\Omega_{bl}^-)$, resp. $u \in W^{1,p}(\Omega_{bl})$, satisfying $u\vert_{\Gamma_{bl}} = 0$, one has 
$$ \| \na u \|_{L^p(\Omega_{bl}^-)}  \le C \|  D u \|_{L^p(\Omega_{bl}^-)}, \quad \mbox{resp.} \: \| \na u \|_{L^p(\Omega_{bl})}  \le C \|  D u \|_{L^p(\Omega_{bl})}. $$
\end{Corollary}
{\em Proof}. One can extend $u$ by $0$ for all $y$ with $-1 < y_2 < \gamma(y_1)$, and apply the previous inequality on $S_{-1}$, resp. $H_{-1}$.

\begin{Proposition}[Rescaled trace and Poincar\'e inequalities]\label{rescaledTracePoincare}
Let $\varphi \in  W^{1,p}(R^\ep)$. We have
	\begin{equation}\label{IQ1}
	\| \varphi \|_{L^{p}(\Sigma)} \leq C \ep^{\frac{1}{p'}} \| \nabla_x \varphi \|_{L^p(R_\ep)} ,
	\end{equation}
	\begin{equation}\label{IQ2}
	\| \varphi \|_{L^p(R_\ep)} \leq C \ep \| \nabla_x \varphi \|_{L^p(R_\ep)} .
	\end{equation}
\end{Proposition}
{\em Proof.} Let $\tilde\varphi(y) = \varphi(\ep y)$, where $y\in S_k = S + (k,-1)$ (a rescaled single cell of rough layer). 
Then $\tilde \varphi \in W^{1,p}(S_k)$ for all $k\in \N$, and $\varphi = 0 $ on $\Gamma$. By the trace theorem and the Poincar\'e 
inequality: for all $p \in [1,\infty )$
	$$\int_{S_k \cap \{y_2 = 0 \}} | \tilde\varphi(\bar{y},0)|^p {\rm\,d}\bar{y} \leq C \int_{S_k} |\nabla_y \tilde \varphi |^p {\rm\,d}y .$$
A change of variables provides
	$$\int_{\ep S_k \cap \{x_2 = 0 \}} | \varphi(\bar{x},0)|^p \ep^{-1} {\rm\,d}\bar{x} 
	\leq C \int_{\ep S_k} \ep^p |\nabla_x \tilde\varphi(x) |^p \ep^{-2} {\rm\,d}x .$$
Summing over $k$ we obtain
	$$\left( \int_{\Sigma} | \varphi(\tilde{x},0) |^p {\rm\, d} \tilde{x} \right)^{\frac{1}{p}} 
	\leq C \ep^{\frac{p-1}{p}} \left( \int_{R^\ep} | \nabla_{x} \varphi (x) |^p \dx  \right)^{\frac{1}{p}}
	$$
and \eqref{IQ1} is proved. The inequality \eqref{IQ2} is proved in the same way, as a consequence of the (one-dimensional) Poincar\'e inequality.  
\qed

%

\end{document}